\documentclass[11pt]{article}
\usepackage{amsmath, amssymb, amscd, epsfig}

\newtheorem{thm}{Theorem}
\newtheorem{prop}[thm]{Proposition}
\newtheorem{lem}[thm]{Lemma}

\newtheorem{cor}[thm]{Corollary}
\newtheorem{rem}[thm]{Remark}
\newtheorem{df}[thm]{Definition}
\newtheorem{ex}[thm]{Example}
\newtheorem{conditions}[thm]{Conditions}

\renewcommand{\phi}{\varphi}
\renewcommand{\epsilon}{\varepsilon}

\newcommand{\BB}{\mathbb}
\newcommand{\g}{\mathfrak}
\newcommand{\pf}{\noindent {\it Proof. }}
\newcommand{\qed}{\nopagebreak $\qquad$ $\square$ \vskip5pt}
\newcommand{\separate}{\vskip5pt}
\newcommand{\den}{\operatorname{Den}}
\newcommand{\talpha}{\widetilde{\alpha}}
\newcommand{\supp}{\operatorname{supp}}
\newcommand{\im}{\operatorname{Im}}
\newcommand{\re}{\operatorname{Re}}
\newcommand{\C}{\operatorname{C}^b_{\BB R-c}}
\newcommand{\D}{\operatorname{D}^b_{\BB R-c}}
\newcommand{\K}{\operatorname{K} (\D (M))}

\begin{document}

\title{\bf Integrals of Equivariant Forms and
a Gauss-Bonnet Theorem for Constructible Sheaves}
\author{Matvei Libine}
\maketitle

\begin{abstract}
The Berline-Vergne integral localization formula for equivariantly
closed forms (\cite{BV}, Theorem 7.11 in \cite{BGV})
is well-known and requires the acting Lie group to be compact.
It is restated here as Theorem \ref{BGV}.
In this article we extend this result to real reductive Lie groups $G_{\BB R}$.
The main result is Theorem \ref{main}.

As an application of this generalization, we prove an analogue of the
Gauss-Bonnet theorem for constructible sheaves (Theorem \ref{bonnet}).
If ${\cal F}$ is a $G_{\BB R}$-equivariant sheaf on a complex
projective manifold $M$, then the Euler characteristic of $M$
with respect to ${\cal F}$
$$
\chi (M,{\cal F}) = \frac 1{(2\pi)^{\dim_{\BB C} M}}
\int_{Ch({\cal F})} \widetilde{\chi_{\g g_{\BB C}}}
$$
as distributions on $\g g_{\BB R}$, where
$Ch({\cal F})$ is the characteristic cycle of ${\cal F}$ and
$\widetilde{\chi_{\g g_{\BB C}}}$ is the Euler form of $M$ extended to
the cotangent space $T^*M$ (independently of ${\cal F}$).
We also consider an analogue of Duistermaat-Heckman measures for
real reductive Lie groups acting on symplectic manifolds.

In \cite{L3} I apply the results of this article to obtain a
Riemann-Roch-Hirzebruch type integral formula for characters of
representations of reductive groups.
\end{abstract}

\noindent {\bf Keywords:}
equivariant forms, Berline-Vergne integral localization formula,
characteristic cycles of sheaves, integral character formula.

\separate

\tableofcontents

\newpage

\begin{section}
{Introduction}  \label{intro}
\end{section}

Equivariant forms were introduced in 1950 by Henri Cartan.
There are many good texts on this subject including \cite{BGV}
and \cite{GS}.

Let $G_{\BB R}$ be a compact Lie group acting on a compact manifold $M$,
let $\g g_{\BB R}$ be the Lie algebra of $G_{\BB R}$,
and let $\alpha(X)$ be an equivariantly closed form on $M$ depending on
$X \in \g g_{\BB R}$.
For $X \in \g g_{\BB R}$, we denote by $M_0(X)$ the set of zeroes of the vector
field on $M$ induced by the infinitesimal action of $X$.
We assume that $M_0(X)$ is discrete.
Then Theorem 7.11 in \cite{BGV} (which we restate here as Theorem \ref{BGV})
says that the integral of $\alpha(X)$ can be expressed as a sum over the
set of zeroes $M_0(X)$ of certain {\em local} quantities of $M$ and $\alpha$:
$$
\int_M \alpha(X) = \sum_{p \in M_0(X)}
\text{local invariant of $M$ and $\alpha$ at $p$}.
$$
This is the essence of the Berline-Vergne integral localization formula for
equivariantly closed differential forms which originally appeared in \cite{BV}.

In this article we extend this result to reductive groups.
So let $G_{\BB R}$ be a real reductive Lie group which may not be compact.
To avoid pathologies we require the action of $G_{\BB R}$ to be complex
algebraic.
On the other hand, for the purpose of interesting applications we would
like to allow integration over homology cycles with non-compact support.
%In order to have a truly new localization formula for non-compact group
%actions which cannot be reduced to the known case of compact group actions
%one must allow integration over cycles with non-compact support.
Then one encounters the following two problems.
First of all, the cycle being infinite, the integral may no longer converge
in the usual sense. We resolve this problem by defining a new
(more relaxed) notion of integral over the cycle in the sense of distribution
on the Lie algebra $\g g_{\BB R}$.
Secondly, some cycles simply may not contain enough points fixed by the
group action for an integral localization formula to make sense.
This is similar to the failure of the Lefschetz fixed point formula for
non-compact manifolds -- some fixed points may run off to infinity.
For this reason we restrict ourselves to the following setting.
Let $G_{\BB R}$ act algebraically on a complex projective manifold $M$,
this action extends naturally to the cotangent space $T^*M$.
Let $\Lambda$ be a conic $G_{\BB R}$-invariant Lagrangian cycle
$\Lambda$ in $T^*M$. We describe a class of differential forms
$\talpha(X)$ on $T^*M$ depending on $X \in \g g_{\BB R}$ and define
$\int_{\Lambda} \talpha (X)$ as a distribution on $\g g_{\BB R}$.
Let $\{x_1,\dots,x_d\}$ be the set of zeroes of the vector field on $M$
induced by the infinitesimal action of $X$.
The main result (Theorem \ref{main}) says that this distribution is
given by integration against a function $F$ on $\g g_{\BB R}$ and
\begin{equation}  \label{F}
F(X) = \sum_{k=1}^d
m_k(X) \cdot
\begin{pmatrix}
\text{the contribution of $x_k$ to the} \\
\text{Berline-Vergne localization formula}
\end{pmatrix},
\end{equation}
where $m_k(X)$ is a certain integer multiplicity which
is exactly the local contribution of $x_k$ to the Lefschetz
fixed point formula, as generalized to sheaf cohomology by
M.~Goresky and R.~MacPherson \cite{GM}.
These multiplicities will be determined in terms of
local cohomology of ${\cal F}$, where ${\cal F}$ is any sheaf
with characteristic cycle $Ch({\cal F}) = \Lambda$.
Existence of such a localization formula was conjectured by W.~Schmid
in \cite{Sch}.

The idea is to observe that the integrand is a
closed form (Lemma \ref{closed}),
to pick a sufficiently regular element $X \in \g g_{\BB R}$
and to deform $\Lambda$ into a simple-looking cycle of the following kind:
$$
m_1(X) T^*_{x_1}M +\dots+ m_d(X) T^*_{x_d}M,
$$
where $m_1(X),\dots,m_d(X)$ are the integer multiplicities from (\ref{F})
%$x_1,\dots,x_d$ are the zeroes of the vector field on $M$
%induced by the infinitesimal action of $X$,
and each cotangent space $T^*_{x_k}M$ is given a certain orientation.
The cycles in question have infinite support, which means one must deform
$\Lambda$ very carefully to ensure that the integral stays unchanged.
The precise result is stated in Proposition \ref{C(X)prop}.

This kind of argument fits very well into the cobordism theory of
spaces equipped with abstract moment maps as described by
V.~Guillemin, V.~Ginzburg and Y.~Karshon in \cite{GGK}.
They would probably call Proposition \ref{C(X)prop}
``the linearization theorem for characteristic cycles.''
Then Theorem \ref{main} essentially becomes ``linearization commutes
with integration.'' Of course, since we work with cycles with possibly
singular support we no longer require that the chains realizing cobordisms
have smooth support.

%While a general cycle may not contain enough fixed points, sometimes one can
%deform such a cycle into a conic Lagrangian cycle for which the integral
%localization formula is known to be true.
%I make an attempt to study such deformations in \cite{L3}.

Then, using this generalized localization formula,
we prove an analogue of the
Gauss-Bonnet theorem for constructible sheaves (Theorem \ref{bonnet}).
If ${\cal F}$ is a $G_{\BB R}$-equivariant sheaf on a complex
projective manifold $M$, then the Euler characteristic of $M$
with respect to ${\cal F}$
$$
\chi (M,{\cal F}) = \frac 1{(2\pi)^{\dim_{\BB C} M}}
\int_{Ch({\cal F})} \widetilde{\chi_{\g g_{\BB C}}}
$$
as distributions on $\g g_{\BB R}$, where
$Ch({\cal F})$ is the characteristic cycle of ${\cal F}$ and
$\widetilde{\chi_{\g g_{\BB C}}}$ is the Euler form of $M$ extended to
the cotangent space $T^*M$ (independently of ${\cal F}$).

In the last section we describe an analogue of Duistermaat-Heckman measures
for real reductive Lie groups acting on symplectic manifolds and give a
formula for the Fourier transforms of these measures similar to the
exact stationary phase approximation formula (Proposition \ref{DH}).

In \cite{L3} I apply the results of this article to obtain a
Riemann-Roch-Hirzebruch type integral formula for characters of
representations of reductive groups.

The proof given here is a significant modification of the
localization argument which appeared in my Ph.D. thesis \cite{L1}.
This thesis provides a geometric proof of an analogue of Kirillov's
character formula for reductive Lie groups.
Article \cite{L2} gives a very accessible introduction to \cite{L1}
and explains key ideas used there by way of examples.

\separate

The following convention will be in force throughout these notes:
whenever $A$ is a subset of $B$,
we will denote the inclusion map $A \hookrightarrow B$ by
$j_{A \hookrightarrow B}$.

\separate

\begin{section}
{The Berline-Vergne Localization Formula}  \label{BVformula}
\end{section}

In this article we use the same notations as in \cite{BGV}.

\separate

Let $M$ be a ${\cal C}^{\infty}$-manifold of dimension $m$ with
an action of a (possibly non-compact) Lie group $G_{\BB R}$,
and let $\g g_{\BB R}$ be the Lie algebra of $G_{\BB R}$.
The group $G_{\BB R}$ acts on ${\cal C}^{\infty}(M)$ by the formula
$(g \cdot \phi) (x) = \phi(g^{-1} x)$.
For $X \in \g g_{\BB R}$, we denote by $X_M$ the vector field on $M$ given by
(notice the minus sign)
$$
(X_M \cdot \phi) (x) =
\frac d{d\epsilon} \phi \bigl( \exp (-\epsilon X)x \bigr) \Bigr|_{\epsilon=0}.
$$
This way
$$
[X,Y]_M = [X_M, Y_M], \qquad \text{for all $X,Y \in \g g_{\BB R}$,}
$$
which would not be true without this choice of signs.

Let $\Omega^*(M)$ denote the (graded) algebra of smooth complex-valued
differential forms on $M$,
and let ${\cal C}^{\infty}(\g g_{\BB R}) \hat \otimes \Omega^*(M)$
denote the algebra of all smooth $\Omega^*(M)$-valued functions on
$\g g_{\BB R}$.
The group $G_{\BB R}$ acts on an element
$\alpha \in {\cal C}^{\infty}(\g g_{\BB R}) \hat \otimes \Omega^*(M)$
by the formula
$$
(g \cdot \alpha)(X) = g \cdot (\alpha ( g^{-1} \cdot X))
\qquad \text{for all $g \in G$ and $X \in \g g_{\BB R}$.}
$$
Let $\Omega^*_{G_{\BB R}}(M) =
({\cal C}^{\infty}(\g g_{\BB R}) \hat \otimes \Omega^*(M))^{G_{\BB R}}$
be the subalgebra of $G_{\BB R}$-invariant elements.
An element $\alpha$ of $\Omega^*_{G_{\BB R}}(M)$ satisfies the relation
$ \alpha(g \cdot X) = g \cdot \alpha(X)$ and is called an
{\em equivariant differential form}.

We define the equivariant exterior differential $d_{\g g_{\BB R}}$ on
${\cal C}^{\infty}(\g g_{\BB R}) \hat \otimes \Omega^*(M)$ by the formula
$$
(d_{\g g_{\BB R}} \alpha) (X) = d(\alpha(X)) - \iota(X_M) (\alpha(X)),
$$
where $d$ denotes the ordinary de Rham differential and
$\iota(X_M)$ denotes contraction by the vector field $X_M$.
This differential $d_{\g g_{\BB R}}$ preserves $\Omega^*_{G_{\BB R}}(M)$, and
$(d_{\g g_{\BB R}})^2 \alpha =0$ for all $\alpha \in \Omega^*_{G_{\BB R}}(M)$.
The elements of $\Omega^*_{G_{\BB R}}(M)$ such that
$d_{\g g_{\BB R}} \alpha =0$ are called {\em equivariantly closed forms}.

\separate

\begin{ex} \label{sigma_R}
{\em Let $T^*M$ be the cotangent bundle of $M$, and let
$\pi: T^*M \twoheadrightarrow M$ denote the projection map.
Let $\sigma_{\BB R}$ denote the canonical symplectic form on $T^*M$.
It is defined, for example, in \cite{KaScha}, Appendix A2.
The action of the Lie group $G_{\BB R}$ on $M$ naturally extends
to $T^*M$. Then we always have a canonical equivariantly closed
form on $T^*M$, namely, $\mu_{\BB R} + \sigma_{\BB R}$.
Here $\mu_{\BB R}: \g g_{\BB R} \to {\cal C}^{\infty}(T^*M)$ is the moment map
defined by:
$$
\mu_{\BB R}(X): \zeta \mapsto -\langle \zeta, X_M \rangle,
\qquad X \in \g g_{\BB R},\: \zeta \in T^*M. \qquad \qquad \square
$$}
\end{ex}

If $\alpha$ is a non-homogeneous equivariant differential form,
$\alpha_{[k]}$ denotes the homogeneous component of degree $k$.
If $M$ is a compact oriented manifold, we can integrate equivariant
differential forms over $M$, obtaining a map
$$
\int_M : \Omega^*_{G_{\BB R}}(M) \to
{\cal C}^{\infty}(\g g_{\BB R})^{G_{\BB R}},
$$
by the formula $(\int_M \alpha) (X) = \int_M \alpha(X)_{[\dim M]}$.

Notice that if $\alpha \in \Omega^*_{G_{\BB R}}(M)$ has top homogeneous
component $\alpha_{[k]}$, then
$(d_{\g g_{\BB R}} \alpha)(X)_{[k+1]}$ is exact;
and if $p \in M$ is a zero of the vector field $X_M$ (i.e. $X_M(p)=0$),
then $(d_{\g g_{\BB R}} \alpha)_{[0]}(p)=0$.
Hence the map $\alpha \mapsto \alpha(X)_{[0]}(p)$ descends to
$\Omega^*_{G_{\BB R}}(M) / \im(d_{\g g_{\BB R}})$.
Similarly, if $M$ is compact, then the map
$\alpha \mapsto \int_M \alpha(X)$ also descends to
$\Omega^*_{G_{\BB R}}(M) / \im(d_{\g g_{\BB R}})$.

Also notice that if $\alpha$ is an equivariantly closed form
whose top homogeneous component has degree $k$, then
$\alpha(X)_{[k]}$ is closed with respect to the ordinary
exterior differential.

\separate

We recall some more notations from \cite{BGV}.
Let $M_0(X)$ be the set of zeroes of the vector field $X_M$.
We state the localization formula in the important
special case where $X_M$ has isolated zeroes.
Here, at each point $p \in M_0(X)$, the Lie action
${\cal X} \mapsto {\cal L}(X_M){\cal X} =[X_M,{\cal X}]$
on the vector fields ${\cal X}$ of $M$
gives rise to a linear transformation $L_p$ on $T_pM$.

If the Lie group $G_{\BB R}$ is compact, then the transformation $L_p$
is invertible and has only imaginary eigenvalues.
Thus the dimension of $M$ is even and there exists an oriented basis
$\{e_1,\dots, e_m\}$ of $T_pM$ such that for $1\le i \le n = m/2$,
$$
L_p e_{2i-1} = \lambda_{p,i} e_{2i}, \quad
L_p e_{2i} = -\lambda_{p,i} e_{2i-1}.
$$
We have $\det(L_p)=\lambda_{p,1}^2 \lambda_{p,2}^2 \dots \lambda_{p,n}^2$,
and it is natural to take the following square root (dependent only on the
orientation of the manifold):
$$
{\det}^{1/2}(L_p) = \lambda_{p,1} \dots \lambda_{p,n}.
$$

\separate

For convenience, we restate Theorem 7.11 from \cite{BGV}.

\begin{thm}  \label{BGV}
Let $G_{\BB R}$ be a compact Lie group with Lie algebra $\g g_{\BB R}$
acting on a compact oriented manifold $M$, and let $\alpha$
be an equivariantly closed differential form on $M$.
Let $X \in \g g_{\BB R}$ be such that the vector field $X_M$ has
only isolated zeroes. Then
$$
\int_M \alpha(X) = (-2\pi)^n
\sum_{p \in M_0(X)} \frac {\alpha(X)_{[0]} (p)}{\det^{1/2} (L_p)},
$$
where $n= \dim(M)/2$, and by $\alpha(X)_{[0]}(p)$, we mean the value of
the function $\alpha(X)_{[0]}$ at the point $p \in M$.
\end{thm}

\separate

\begin{section}
{A Brief Introduction to Characteristic Cycles of Sheaves}
\end{section}

Characteristic cycles were introduced by M.~Kashiwara and their
definition can be found in \cite{KaScha}.
A comprehensive treatment of characteristic cycles can be found in \cite{Schu}.
On the other hand, W.~Schmid and K.~Vilonen give a geometric way to understand
characteristic cycles in \cite{SchV1} which we follow here.
In this section we briefly outline the defining properties of
characteristic cycles which are analogous to
Eilenberg-Steenrod homology axioms for homology of topological spaces.
Let ${\cal F}$ be a sheaf on a manifold $M$.
The characteristic cycle $Ch({\cal F})$ is a conic Lagrangian Borel-Moore
homology cycle lying inside the cotangent space $T^*M$.
If the sheaf ${\cal F}$ happens to be perverse,
the characteristic cycle of ${\cal F}$ equals the characteristic cycle of the
holonomic ${\cal D}$-module corresponding to ${\cal F}$ via the
Riemann-Hilbert correspondence.

In this section only we assume that $M$ is an oriented smooth real
semi-algebraic manifold which need not be compact.
See, for instance, \cite{DM} for the notion of semi-algebraic sets.
(In the next section we will further require $M$ to be a smooth complex
projective variety).
%Recall that a sheaf ${\cal F}$ on $M$ is {\em locally constant} if there is
%an open covering $M = \cup_j U_j$ such that for each $j$ the restricted sheaf
%${\cal F} |_{U_j}$ is constant.
Now let ${\cal F}$ be a bounded complex of sheaves on $M$.
We say that ${\cal F}$ has {\em $\BB R$-constructible cohomology}
if there exists a locally finite covering $M = \cup_{j \in J} M_j$
by semi-algebraic subsets such that for all $k \in \BB Z$ and all $j \in J$,
the restricted cohomology sheaves $H^k({\cal F}) |_{M_j}$ are
constant of finite rank.

Let $\C (M)$ denote the category of bounded complexes of sheaves on $M$
with $\BB R$-constructible cohomology, and let $\D (M)$ denote the
bounded derived category of sheaves on $M$ with $\BB R$-constructible
cohomology. From now on ${\cal F}$ denotes an element in $\C (M)$ or $\D (M)$.
The characteristic cycle
$Ch({\cal F})$ associated to ${\cal F}$ is a Borel-Moore homology cycle
(possibly with infinite support) in the cotangent space $T^*M$
of dimension $\dim_{\BB R}M$.
The cycle $Ch({\cal F})$ has the following properties: it is
conic, i.e. invariant under the scaling action of positive reals
$\BB R^{>0}$ on $T^*M$
(but not necessarily under the action of $\BB R^{\times}$),
and its support $|Ch({\cal F})|$ is Lagrangian.
More precisely, there exists a Whitney
stratification ${\cal S}$ of $M$ by semi-algebraic sets
such that the cohomology of ${\cal F}$
is constructible with respect to ${\cal S}$.
This means that for all $k \in \BB Z$ and all $S \in {\cal S}$,
the cohomology sheaves restricted to the strata
$H^k({\cal F}) |_S$ are (locally) constant of finite rank.
Then the support of $Ch({\cal F})$ lies in the union of conormal spaces:
$$
|Ch({\cal F})| \subset \bigcup_{S \in {\cal S}} T^*_SM.
$$

Let ${\cal L}^+(M)$ denote the abelian group of Borel-Moore homology cycles
(with coefficients in $\BB Z$) in the cotangent space $T^*M$ of dimension
$\dim_{\BB R}M$ which are conic (i.e. invariant under the scaling action of
$\BB R^{>0}$ on $T^*M$) and whose support lies in
$\cup_{S \in {\cal S}} T^*_SM$ for some locally finite semi-algebraic
Whitney stratification ${\cal S}$ of $M$.

\begin{ex}
{\em
Let $N \subset M$ be a closed semi-algebraic submanifold,
$j: N \hookrightarrow M$ the inclusion map, and let
$\BB C_N$ be the constant sheaf on $N$ of rank 1, then
$Ch(j_* \BB C_N)$ is the conormal space $T^*_NM$ equipped
with a certain orientation.

To specify this orientation, pick any point $p \in N$ and choose a
positively oriented system of coordinates $(x_1,\dots,x_{\dim M})$
on $M$ around $p$ such that $N$ is locally given by the equations
$x_{\dim N +1} = \dots = x_{\dim M} = 0$.
Let $(\xi_1,\dots,\xi_{\dim M})$ be the fiber coordinates
dual to the frame $dx_1,\dots,dx_{\dim M}$.
Then near points lying in the cotangent space $T^*_pM$,
$T^*_NM$ is given by the equations
$x_{\dim N +1} = \dots = x_{\dim M} = \xi_1 = \dots = \xi_{\dim N} =0$
and the functions
$(x_1,\dots,x_{\dim N}, \xi_{\dim N +1},\dots,\xi_{\dim M})$
form a coordinate system on $T^*_NM$.
Finally, $Ch(j_* \BB C_N)$ is the conormal space $T^*_NM$ with orientation
equal $(-1)^{\dim M - \dim N}$ times the orientation given by coordinates
$(x_1,\dots,x_{\dim N}, \xi_{\dim N +1},\dots,\xi_{\dim M})$,
and this orientation does not depend on the choices made.
\qed
}\end{ex}

Following W.~Schmid and K.~Vilonen we introduce the notions of families of
cycles and their limits. Let $\tilde M$ be an ambient manifold which we later
take equal $T^*M$, and let $I=(0,b)$ be an open interval.

\begin{df}  \label{family}
A {\em family of $k$-cycles} in $\tilde M$ parametrized by $I$ is a
$(k+1)$-cycle $C_I$ in $I \times \tilde M$ with the following property:
for each $s \in I$, there exists a Whitney stratification of $|C_I|$,
such that the ``slice'' $|C_I| \cap (\{s\} \times \tilde M)$ is
a Whitney stratified subset of $|C_I|$ of dimension at most $k$.
\end{df}

For each $s \in I$, we identify $\tilde M$ with $\{s\} \times \tilde M$
and we have a {\em specialization} map $C_I \mapsto C_s$, where
$C_s$ is a $k$-cycle in $\tilde M \simeq \{s\} \times \tilde M$
with support lying in $|C_I| \cap (\{s\} \times \tilde M)$.
The precise definition can be found in \cite{SchV1},
but we skip it because this notion is quite intuitive and in this article
all families of cycles will be defined explicitly through the specializations
$C_s$. Note that if the dimension of $|C_I| \cap (\{s\} \times \tilde M)$
is strictly less than $k$, then $C_s=0$.

Next we define the limit of a family of cycles as the parameter $s \to 0^+$.
Recall that $I$ is an open interval $(0,b)$ and set $J=[0,b)$.
We consider a family of $k$-cycles $C_I$ in $\tilde M$ subject to an
additional condition: the closure $\overline{|C_I|}$ in $J \times \tilde M$
admits a Whitney stratification such that
$\overline{|C_I|} \cap (\{0\} \times \tilde M)$
is a stratified subset of $\overline{|C_I|}$.
Note that $\overline{|C_I|}$ is a subset of $J \times \tilde M$
and the latter is a manifold with boundary, so to make sense out of its
Whitney stratification we embed $J \times \tilde M$ into
$\BB R \times \tilde M$.
Then it follows that $\overline{|C_I|} \cap (\{0\} \times \tilde M)$
has dimension at most $k$. The $(k+1)$-cycle $C_I$ in $I \times \tilde M$
can be regarded as a $(k+1)$-chain in $J \times \tilde M$, the boundary of
this chain $\partial C_I$ is necessarily supported in $\{0\} \times \tilde M$.
Since $\{0\} \times \tilde M \simeq \tilde M$, we regard $\partial C_I$
as a cycle in $\tilde M$ and define
$$
\lim_{s \to 0^+} C_s = -\partial C_I.
$$
The negative sign appears for orientation reasons and ensures that the
formal definition of a limit agrees with geometric intuition behind it.

\begin{prop}[Proposition 3.25 in \cite{SchV1}]  \label{deform}
For all $t \in I$,
$$
C_t - (\lim_{s \to 0^+} C_s) = \partial C_{(0,t)},
$$
where $C_{(0,t)}$ denotes the restriction of $C_I$ to $(0,t) \times \tilde M$.
\end{prop}

Let $U$ be an open semi-algebraic subset in $M$.
We are going to define two pushforward maps of cycles
${\cal L}^+(U) \to {\cal L}^+(M)$.
By a semi-algebraic function on $M$ we mean a function whose
graph is a semi-algebraic subset of $M \times \BB R$.
Then one can find a semi-algebraic function $f$ on $M$ of class ${\cal C}^2$
such that $f$ is strictly positive on $U$ and the boundary
$\partial U$ is {\rm }precisely the zero set of $f$
(Proposition 4.22 in \cite{DM}).
Let $j$ denote the inclusion map $U \hookrightarrow M$, and let
$\Lambda \in {\cal L}^+(U)$ be a conic Lagrangian cycle in $T^*U$.
For each $s>0$, we regard $s \frac {df}f$ as a section in $T^*U$,
it induces two mutually inverse homeomorphisms of $T^*U$ defined fiberwise by:
$$
\tau_+: (\zeta, x) \mapsto \Bigl( \zeta + s \frac {df}f(x), x \Bigr),
\qquad
\tau_-: (\zeta, x) \mapsto \Bigl( \zeta - s \frac {df}f(x), x \Bigr).
$$
We set
$$
\Lambda + s \frac {df}f = (\tau_+)_* (\Lambda),
\qquad
\Lambda - s \frac {df}f = (\tau_-)_* (\Lambda).
$$

\begin{prop}
Under the above hypotheses,
the cycles $\Lambda \pm s \frac {df}f$ in $T^*U$, regarded as chains
in $T^*M$, have no boundary, they form two families of cycles in $T^*M$
parametrized by $(0,\infty)$, and the limits
$$
\lim_{s \to 0^+} \Bigl( \Lambda + s \frac {df}f \Bigr),
\qquad
\lim_{s \to 0^+} \Bigl( \Lambda - s \frac {df}f \Bigr).
$$
do not depend on the choice of a semi-algebraic function $f$ on $M$ of class
${\cal C}^2$ such that $f>0$ on $U$ and the zero set
$\{f = 0 \} = \partial U$.
\end{prop}

This proposition can be extracted from Section 4 of \cite{SchV1}.
The growth of $\frac {df}f$ near the boundary of $U$ ensures that
$\Lambda \pm s \frac {df}f$ are cycles in $T^*M$ as opposed to chains
with boundary.
The proposition implies that the following two maps are well-defined:
$$
j_*: {\cal L}^+(U) \to {\cal L}^+(M), \qquad \Lambda \mapsto 
\lim_{s \to 0^+} \Bigl( \Lambda + s \frac {df}f \Bigr)
$$
and
$$
j_!: {\cal L}^+(U) \to {\cal L}^+(M), \qquad \Lambda \mapsto 
\lim_{s \to 0^+} \Bigl( \Lambda - s \frac {df}f \Bigr).
$$

\begin{ex}
{\em
Let $M = \BB R$ with its standard orientation, and let $U = (0,\infty)$.
Take $\Lambda \in {\cal L}^+(0,\infty)$ equal the zero section of $T^*U$
oriented the same way $U$ is. This $\Lambda$ is the characteristic cycle
of $\BB C_{(0,\infty)}$ -- the constant sheaf on $(0,\infty)$ of rank 1.
Note that $\Lambda$, regarded as a chain in $T^* \BB R$, has non-zero boundary.
We can take the defining function of $(0,\infty)$ to be $f(x)=x$,
where $x$ is the standard coordinate on $\BB R$.
Then $\Lambda + s \frac{df}{f} = s \frac{dx}x$ is a piece of hyperbola
$\xi = \frac sx$, $x >0$, oriented so that the direction of increasing $x$
is positive. Evidently, these cycles have no boundary in $T^* \BB R$.
As $s \to 0^+$, these cycles tend to $j_*(\Lambda)$ which has support
$\{\xi > 0 \} \cup \{ x>0 \}$ and oriented along decreasing $\xi$ and
increasing $x$.
Similarly, $\Lambda - s \frac{df}{f} = -s \frac{dx}x$ is a piece of hyperbola
$\xi = -\frac sx$, $x >0$, oriented along increasing $x$.
As $s \to 0^+$, these cycles tend to $j_!(\Lambda)$ which has support
$\{\xi < 0 \} \cup \{ x>0 \}$ and oriented along increasing $\xi$ and
increasing $x$.

We will see in a moment that the limit cycles
$j_*(\Lambda)$ and $j_!(\Lambda)$ are the characteristic
cycles of $Rj_*(\BB C_{(0,\infty)})$ and $Rj_!(\BB C_{(0,\infty)})$
respectively.
\qed
}\end{ex}

\begin{df}
The characteristic cycle is a map
$$
Ch: \: \C (M) \to {\cal L}^+(M)
$$
which is uniquely determined by the following properties:
\begin{enumerate}
\item
\underline{Normalization}:
Let $\BB C_M$ be the constant sheaf on $M$ of rank 1, then
$$
Ch(\BB C_M) = [M] =
\begin{matrix}
\text{zero section of $T^*M$ oriented} \\
\text{by the fixed orientation of $M$}
\end{matrix} \,;
$$

\item
\underline{Additivity}:
The map $Ch$ descends to $\D (M)$ -- the bounded derived category of
sheaves on $M$ with $\BB R$-constructible cohomology -- and is additive on
distinguished triangles of $\D (M)$:
$$
Ch({\cal F}) = Ch({\cal F}') + Ch({\cal F}'')
$$
whenever there is a distinguished triangle
$$
\begin{CD}
{\cal F}' @>>> {\cal F} @>>> {\cal F}'' @>>{+1}>
\end{CD}
$$
in $\D(M)$;

\item
\underline{$Ch$ Is Local}:
For any open semi-algebraic subset $U \subset M$ the following diagram
commutes:
$$
\begin{CD}
\C (M)   @>{Ch}>>   {\cal L}^+(M) \\
@VVV      @VVV \\
\C (U)    @>{Ch}>> {\cal L}^+(U),
\end{CD}
$$
where the left vertical arrow is the restriction map of complexes of sheaves
and the right vertical arrow is the restriction map of cycles with infinite
support from $T^*M$ to its open subset $T^*U$;

\item
\underline{Open Embedding}:
For any open semi-algebraic subset $U \subset M$ the following diagram
commutes:
$$
\begin{CD}
\C (U)   @>{Ch}>>   {\cal L}^+(U) \\
@V{R_*j}VV      @VV{j_*}V \\
\C (M)    @>{Ch}>> {\cal L}^+(M).
\end{CD}
$$
\end{enumerate}
\end{df}

As was explained in \cite{SchV1}, these properties uniquely determine the
characteristic cycle map $Ch: \: \C (M) \to {\cal L}^+(M)$, however, from
this point of view proving its existence becomes a highly non-trivial matter.
Below we state more properties of characteristic cycles, starting with a
stronger open embedding property.

\begin{thm}[Open Embedding Theorem 4.2 in \cite{SchV1}]  \label{oet}
Let $U$ be an open semi-algebraic subset in $M$, and let $f$
be semi-algebraic function on $M$ of class ${\cal C}^2$
such that $f>0$ on $U$ and the boundary $\partial U$ is precisely
the zero set of $f$.
Let ${\cal F} \in \C (U)$ be a bounded complex of sheaves on $U$ with
${\BB R}$-constructible cohomology. Then
\begin{align*}
Ch(Rj_* {\cal F}) &= \lim_{s \to 0^+} \Bigl( Ch({\cal F}) +s \frac {df}f \Bigr)
=j_* (Ch({\cal F})),  \\
Ch(Rj_! {\cal F}) &= \lim_{s \to 0^+} \Bigl( Ch({\cal F}) -s \frac {df}f \Bigr)
=j_! (Ch({\cal F})).
\end{align*}
\end{thm}

The Open Embedding Theorem not only provides a means of computing the
characteristic cycles of $Rj_* {\cal F}$ and $Rj_! {\cal F}$, but also
a way of deforming them.
The following observation will play a crucial role in
Section \ref{deformationsection}.
It immediately follows from the Open Embedding Theorem and Proposition
\ref{deform}.

\begin{cor}
The following pairs of cycles are homologous:
$$
Ch(Rj_* {\cal F}) \quad \sim \quad Ch({\cal F}) + \frac {df}f
\qquad \text{and} \qquad
Ch(Rj_! {\cal F}) \quad \sim \quad Ch({\cal F}) - \frac {df}f.
$$
Moreover the chains realizing these homology relations can be chosen to lie
inside the sets
$$
\overline {\bigcup_{0 \le s \le 1} \Bigl(|Ch({\cal F})|+ s \frac {df}f} \Bigr)
\qquad \text{and} \qquad
\overline {\bigcup_{0 \le s \le 1} \Bigl(|Ch({\cal F})|- s \frac {df}f} \Bigr)
$$
respectively.
\end{cor}

Let $\K$ denote the Grothendieck group of $\D(M)$, i.e. the abelian group
generated by the objects of $\D(M)$ with one relation
${\cal F} = {\cal F}'+ {\cal F}''$ for each
distinguished triangle
$$
\begin{CD}
{\cal F}' @>>> {\cal F} @>>> {\cal F}'' @>>{+1}>
\end{CD}
$$
in $\D(M)$.
The additivity property of characteristic cycles implies that $Ch$ descends to
a homomorphism
$$
Ch: \: \K \to {\cal L}^+(M).
$$
M.~Kashiwara and P.~Schapira (Theorem 9.7.10 in \cite{KaScha}) show that 
this homomorphism is in fact an isomorphism of abelian groups.
In particular, every conic Lagrangian cycle in $T^*M$ can be realized
as the characteristic cycle of some ${\cal F} \in \C (M)$.
If a group $G_{\BB R}$ acts on $M$ semi-algebraically and
${\cal F} \in \C (M)$ is $G_{\BB R}$-equivariant
(see \cite{KaSchm} for the definition),
then $Ch({\cal F})$ is $G_{\BB R}$-invariant. Furthermore,
$$
\langle \mu_{\BB R}(\zeta), X \rangle
= - \langle \zeta, X_M \rangle = 0
\qquad \text{for all $\zeta \in |Ch({\cal F})|$, $X \in \g g_{\BB R}$},
$$
where the vector field $X_M$ and the real moment map $\mu_{\BB R}$ are defined
in Section \ref{BVformula}.
Conversely, every $G_{\BB R}$-invariant cycle $\Lambda \in {\cal L}^+(M)$
can be realized as the characteristic cycle of some
$G_{\BB R}$-equivariant ${\cal F} \in \C (M)$.

\begin{thm} [Hopf Index Theorem (Corollary 9.5.2 in \cite{KaScha})] \label{hit}
Suppose that a complex of sheaves ${\cal F} \in \C (M)$ has compact support,
then the Euler characteristic of $M$ with respect to ${\cal F}$
$$
\chi(M,{\cal F}) = \#([M] \cap Ch({\cal F})),
$$
where the right hand side denotes the intersection number between the
cycles $[M]$ and $Ch({\cal F})$.
\end{thm}

For $k \in \BB Z$, let ${\cal F}[k]$ denote the complex ${\cal F}$ with
degrees shifted by $k$, then $Ch({\cal F}[k]) = (-1)^k Ch({\cal F})$.
Let $\BB D_M : \D (M) \to \D (M)$ denote the Verdier duality operator,
then $Ch(\BB D_M({\cal F})) = Ch({\cal F})^a$, where
$a : T^*M \to T^*M$ is the antipodal map $\zeta \mapsto -\zeta$ and
$Ch({\cal F})^a$ denotes the image of $Ch({\cal F})$ under this map.

If $f: M \to N$ is a proper map of real semi-algebraic manifolds,
there is an explicit description of the effect on characteristic cycles
by the pushforward map $Rf_*: \D (M) \to \D (N)$.
Similarly, if $f: M \to N$ is a map of real semi-algebraic manifolds,
${\cal G} \in \D(N)$ and $f$ is ``normally non-singular with respect to
${\cal G}$'' (a transversality condition on the induced map
$df: T^*N \to T^*M$ and the stratification ${\cal S}$ of $M$ making
the cohomology of ${\cal G}$ constant),
there is an explicit description of $Ch(Rf^*({\cal G}))$.
(See, for instance, \cite{KaScha}, \cite{SchV1}).

\separate

\begin{section}
{Statement of the Main Result}
\end{section}

Let $G_{\BB C}$ be a connected complex algebraic reductive group
which is defined over $\BB R$, and let $G_{\BB R}$ be a subgroup of
$G_{\BB C}$ lying between the group of real points $G_{\BB C}(\BB R)$
and the identity component $G_{\BB C}(\BB R)^0$.
We regard $G_{\BB R}$ as a real reductive Lie group.
Let $\g g_{\BB C}$ and $\g g_{\BB R}$ be their respective Lie algebras.
We pick another subgroup $U_{\BB R}$ of $G_{\BB C}$ such that, letting
$\g u_{\BB R}$ be the Lie algebra of $U_{\BB R}$, we have an isomorphism
$\g u_{\BB R} \otimes_{\BB R} \BB C \simeq \g g_{\BB C}$.
In the applications we have in mind we will choose $U_{\BB R}$ to be a compact
real form of $G_{\BB C}$ (i.e. a maximal compact subgroup of $G_{\BB C}$),
but we do not require $U_{\BB R}$ to be compact for now;
for instance, $U_{\BB R}$ may equal $G_{\BB R}$.
Let $M$ be a smooth complex projective variety of dimension $n$
with a complex algebraic $G_{\BB C}$-action on it.
We denote by $\Omega^{(p,q)}(M)$ the space of complex-valued differential
forms of type $(p,q)$ on $M$.

Let $T^*M$ be the holomorphic cotangent bundle of $M$, and let
$\pi: T^*M \twoheadrightarrow M$ denote the projection map.
Let $\sigma$ denote the canonical complex algebraic holomorphic symplectic
form on $T^*M$ defined similarly to the form $\sigma_{\BB R}$ from
Example \ref{sigma_R}.
The action of the Lie group $G_{\BB C}$ on $M$ naturally extends
to $T^*M$. Then we always have a canonical equivariantly closed
form on $T^*M$, namely, $\mu + \sigma$.
Here $\mu: \g g_{\BB C} \to {\cal C}^{\infty}(T^*M)$ is the moment map
defined by:
\begin{equation}  \label{moment}
\mu(X): \zeta \mapsto - \langle \zeta, X_M \rangle,
\qquad X \in \g g_{\BB C},\: \zeta \in T^*M.
\end{equation}

\begin{rem}  \label{T*M}
If $M$ is a complex manifold and
$M^{\BB R}$ is the underlying real analytic manifold,
then there are at least two different but equally natural ways 
to identify the real cotangent bundle $T^*(M^{\BB R})$ with the
holomorphic cotangent bundle $T^*M$ of the complex manifold $M$.
We use the convention (11.1.2) of \cite{KaScha}, Chapter XI;
the same convention is used in \cite{L1}, \cite{L2} and \cite{SchV2}.
Under this convention, if
$\sigma_{\BB R}$ is the canonical real symplectic form on $T^*M^{\BB R}$
described in Example \ref{sigma_R} and
$\sigma$ is the canonical complex symplectic form on $T^*M$,
then $\sigma_{\BB R}$ gets identified with $2 \re \sigma$.
(And $\mu_{\BB R} = 2\re \mu$.)
% dz <--> 2dx,  -idz <--> 2dy
\end{rem}

In this article we consider integrals over Borel-Moore homology cycles
$\Lambda$ in $T^*M$ (with coefficients in $\BB Z$)
which satisfy the following three properties:
\begin{itemize}
\item
$\Lambda$ is {\em real-Lagrangian}, i.e.
$\dim_{\BB R} \Lambda = \dim_{\BB R} M$
and there exists a locally finite semi-algebraic Whitney stratification
${\cal S}$ of $M^{\BB R}$ such that, regarding $\Lambda$ as a cycle in
$T^*(M^{\BB R})$ via the identification with $T^*M$, the support
of $\Lambda$ lies in $\cup_{S \in {\cal S}} T^*_SM$;
\item
$\Lambda$ is {\em conic}, i.e. invariant under the scaling action of
positive reals $\BB R^{>0}$ on $T^*M$ (but not necessarily under
the actions of $\BB C^{\times}$ or $\BB R^{\times}$);
\item
$\Lambda$ is $G_{\BB R}$-invariant.
\end{itemize}
We denote the abelian group of such cycles by ${\cal L}^+_{G_{\BB R}} (M)$.
Note that the Lagrangian condition together with
$G_{\BB R}$-equivariance imply $\re \sigma |_{\Lambda} \equiv 0$
and $\mu(|\Lambda|) \subset i \g g_{\BB R}^*$.
As was mentioned earlier, given any $\Lambda \in {\cal L}^+_{G_{\BB R}} (M)$,
there exists a $G_{\BB R}$-equivariant complex of sheaves
${\cal F} \in \C (M)$ such that $\Lambda = Ch({\cal F})$.
The reason for restricting ourselves to the conic Lagrangian cycles in $T^*M$
was explained in Section \ref{intro}.

\begin{ex}
{\em
Consider $G_{\BB R} = GL(l,\BB R) \subset GL(l, \BB C) = G_{\BB C}$
acting naturally on a complex Grassmanian $Gr_{\BB C}(k,l)$.
Let $N$ be the real Grassmanian $Gr_{\BB R}(k,l) \subset Gr_{\BB C}(k,l)$ and
$\Lambda = T^*_{Gr_{\BB R}(k,l)} Gr_{\BB C}(k,l)$
equipped with some orientation.
\qed
}\end{ex}

\begin{conditions}  \label{conditions}  {\em
We consider forms $\alpha : \g g_{\BB C} \to \Omega^*(M)$
which satisfy the following three conditions:
\begin{enumerate}
\item
The assignment $X \mapsto \alpha(X) \in \Omega^*(M)$
depends holomorphically on $X \in \g g_{\BB C}$;
\item
For each $k \in \BB N$ and each $X \in \g g_{\BB C}$,
\begin{equation}  \label{condition2}
\alpha(X)_{[2k]} \in
\bigoplus_{\begin{matrix} p+q=2k \\ p \ge q \end{matrix}} \Omega^{(p,q)}(M);
\end{equation}
\item
The restriction of $\alpha$ to $\g u_{\BB R}$ is an equivariantly closed
form with respect to $U_{\BB R}$.
\end{enumerate}
}\end{conditions}

\begin{ex}  {\em
A $U_{\BB R}$-equivariant characteristic form
$\alpha: \g u_{\BB R} \to \Omega^*(M)$
(defined in Section 7.1 of \cite{BGV}) satisfies the third condition.
Since it depends on $X \in \g u_{\BB R}$ polynomially,
$\alpha$ extends uniquely to a map $\alpha: \g g_{\BB C} \to \Omega^*(M)$
so that the first condition is satisfied.
Finally, for each $X \in \g g_{\BB C}$,
$$
\alpha(X) \in \bigoplus_k \Omega^{(k,k)}(M),
$$
so that the second condition is satisfied too.
This is the most important class of forms satisfying
Conditions \ref{conditions}.
\qed
}\end{ex}

We regard $M$ as a submanifold of $T^*M$ via the zero section inclusion.
We consider the form
%\begin{equation}  \label{talpha}
$$
\talpha(X) = e^{\mu(X) + \sigma} \wedge \pi^* \bigl( \alpha(X) \bigr),
\qquad X \in \g g_{\BB C}.
$$
%\end{equation}
The restriction of $\talpha(X)$ to $M$ is just
$\alpha(X)$. We will see later that, in a way, $\talpha$ is the most
natural equivariant extension of $\alpha$ to $T^*M$.
To avoid cumbersome notations, we denote the image of an element
$\beta \in \Omega^*(M)$ under the inclusion
$\pi^*: \Omega^*(M) \hookrightarrow \Omega^*(T^*M)$ by
$\beta$ as well instead of $\pi^*(\beta)$. Thus
$$
\talpha(X) = e^{\mu(X) + \sigma} \wedge \alpha(X).
$$

%\begin{rem}
%The idea to consider such forms $\talpha: \g g_{\BB C} \to \Omega^*(T^*M)$
%came from the integral character formula due to W.~Schmid and K.~Vilonen
%(\cite{SchV2}). In that character formula the integrand had this shape.
%It was the pullback of a naturally defined form on a complex coadjoint
%orbit to the cotangent space of the flag manifold of $\g g_{\BB C}$
%via the ``twisted moment map.'' 
%\end{rem}

\separate

Recall that $n=\dim_{\BB C} M$, so that the cycle
$\Lambda \in {\cal L}^+_{G_{\BB R}} (M)$ has dimension $2n$.

\begin{lem}  \label{closed}
For each $X \in \g g_{\BB C}$, the form $\talpha(X)_{[2n]}$ is closed.
\end{lem}

\pf
First we show that $\talpha(X)_{[2n+2]}=0$.
Indeed,
$$
\talpha(X)_{[2n+2]}= e^{\mu(X)} \sum_{k=0}^{k=n+1}
\frac 1{(n-k+1)!} \sigma^{n-k+1} \wedge \alpha(X)_{[2k]},
$$
so it suffices to show that each term
$\sigma^{n-k+1} \wedge \alpha(X)_{[2k]} =0$.
But this follows from (\ref{condition2}) and an observation
$$
\sigma^{n-k+1} \wedge \Omega^{(p,q)}(M) =0
\qquad \text{for $p \ge k$.}
$$

The restriction of $\talpha$ to $\g u_{\BB R}$ is equivariantly closed
with respect to the action of $U_{\BB R}$ on $T^*M$ for the reason that
it is ``assembled'' from $U_{\BB R}$-equivariantly closed forms.
Hence $\talpha(X)_{[2n]}$ is closed for all $X \in \g u_{\BB R}$.
But since $d\talpha(X)$ depends on $X \in \g g_{\BB C}$ holomorphically, 
$d\talpha(X)_{[2n]} =0$ for all $X \in \g g_{\BB C}$.
\qed

\separate

If $\phi$ is a smooth compactly supported differential form on
$\g g_{\BB R}$ of top degree, then we define its Fourier transform
as in \cite{L1}, \cite{L2} and \cite{SchV2}:
\begin{equation}  \label{ftransform}
\hat \phi(\xi) = \int_{\g g_{\BB R}} e^{\langle X, \xi \rangle} \phi(X),
\qquad X \in \g g_{\BB R},\,\xi \in \g g_{\BB C}^*,
\end{equation}
without the customary factor of $i=\sqrt{-1}$ in the exponent.

Similarly we define
$\widehat{\phi\alpha} : \g g_{\BB C}^* \to \Omega^*(M)$:
$$
\widehat{\phi\alpha} (\xi) = \int_{\g g_{\BB R}} 
e^{\langle X, \xi \rangle} \phi(X) \wedge \alpha(X),
\qquad X \in \g g_{\BB R},\: \xi \in \g g_{\BB C}^*,
$$
where $\phi(X) \wedge \alpha(X)$ is a form on $\g g_{\BB R} \times M$.
For each $\xi \in \g g_{\BB C}^*$, the form $\widehat{\phi\alpha} (\xi)$
belongs to $\Omega^*(M)$ and decays rapidly as $\xi \to \infty$,
$\xi \in i \g g_{\BB R}^*$.

We can regard the moment map (\ref{moment}) as a map
$\mu: T^*M \to \g g_{\BB C}^*$ via
\begin{equation}  \label{mu}
\mu(\zeta): X \mapsto -\langle \zeta, X_M \rangle,
\qquad X \in \g g_{\BB C},\: \zeta \in T^*M.
\end{equation}
Abusing notation we denote by $\mu^* (\widehat{\phi\alpha}) \in \Omega^*(T^*M)$
the pullback of $\widehat{\phi\alpha} \in \Omega^*(\g g_{\BB C}^* \times M)$
via the composition of
$$
\begin{matrix}
T^*M \hookrightarrow T^*M \times M \\ \zeta \mapsto (\zeta, \pi(\zeta))
\end{matrix}
\qquad \text{and} \qquad
\begin{matrix}
T^*M \times M \to \g g_{\BB C}^* \times M \\
(\zeta, x) \mapsto (\mu(\zeta),x).
\end{matrix}
$$
Then
$$
\mu^* (\widehat{\phi\alpha}) =
\int_{\g g_{\BB R}} 
e^{\langle X, \mu(\zeta) \rangle} \phi(X) \wedge \alpha(X),
\qquad \zeta \in T^*M,\: X \in \g g_{\BB R}.
$$

We will be studying integrals of the kind
\begin{equation}  \label{int}
\int_{\Lambda} \mu^* (\widehat{\phi\alpha}) \wedge e^{\sigma} =
\int_{\Lambda} \Bigl( \int_{\g g_{\BB R}} 
e^{\langle X, \mu(\zeta) \rangle + \sigma}
\wedge \phi(X) \wedge \alpha(X) \Bigr)
= \int_{\Lambda} \Bigl( \int_{\g g_{\BB R}} \talpha \wedge \phi(X) \Bigr).
\end{equation}
Of course, the cycle $\Lambda$ being infinite it is not clear at all
whether this integral converges.
We denote by
$$
\supp(\sigma|_{\Lambda})
$$
the closure in $T^*M$ of the set of smooth points of the support
$|\Lambda|$ where $\sigma|_{|\Lambda|} \ne 0$.

\begin{lem}
If the moment map $\mu: T^*M \to \g g_{\BB C}^*$
is proper on the set $\supp(\sigma|_{\Lambda})$,
then the integral (\ref{int}) converges.
In particular, the integral (\ref{int}) converges
if the moment map $\mu$ is proper on the support $|\Lambda|$.
\end{lem}

\pf
Note that $M$ is compact, so the only unbounded directions of $\Lambda$
are those along the fibers of $T^*M \twoheadrightarrow M$.
We fix any norm $\|.\|_{\g g_{\BB C}^*}$ on $\g g_{\BB C}^*$.
For $R>0$ we denote by $B_R$ the open ball of radius $R$ in $\g g_{\BB C}^*$:
\begin{equation}  \label{B_R}
B_R = \{ \xi \in \g g_{\BB C}^*;\: \|\xi\|_{\g g_{\BB C}^*} <R \}
\end{equation}
and by $\overline{B_R}$ its closure in $\g g_{\BB C}^*$.
We already know that
$$
\widehat{\phi\alpha} (\xi) = \int_{\g g_{\BB R}} 
e^{\langle X, \xi \rangle} \phi(X) \wedge \alpha(X)
$$
decays rapidly as $\|\xi\|_{\g g_{\BB C}^*} \to \infty$,
$\xi \in i \g g_{\BB R}^*$.

Since the cycle $\Lambda$ is $G_{\BB R}$-invariant,
$\mu(|\Lambda|) \subset i \g g_{\BB R}^*$.
On the other hand, $\mu$ being proper on $\supp(\sigma|_{\Lambda})$
implies that the set
$\supp(\sigma|_{\Lambda})  \cap \mu^{-1} (\overline{B_R})$ is compact.
Since the cycle $\Lambda$ is conic along the fibers of
$T^*M \twoheadrightarrow M$ and the integrand decays rapidly
on $\supp(\sigma|_{\Lambda})$ along those fibers as $R \to \infty$,
it is clear that the limit
$$
\lim_{R \to \infty}
\int_{\Lambda \cap \mu^{-1}(B_R)}
\mu^* (\widehat{\phi\alpha}) \wedge e^{\sigma}
= \lim_{R \to \infty}
\int_{\Lambda \cap (M \cup \supp(\sigma|_{\Lambda}))
\cap \mu^{-1}(B_R)} \mu^* (\widehat{\phi\alpha}) \wedge e^{\sigma}
$$
is finite.
\qed

\begin{ex}  {\em
The condition of the lemma is automatically satisfied
if the support $|\Lambda|=M$ (which happens when $\Lambda = Ch(\BB C_M)$,
where $\BB C_M$ is the constant sheaf on $M$ of rank 1).

This condition is also satisfied when $M$ is a
homogeneous space $G_{\BB C} /P_{\BB C}$, where
$P_{\BB C} \subset G_{\BB C}$ is a parabolic subgroup, and
$\Lambda \in {\cal L}^+_{G_{\BB R}} (M)$ is any cycle at all.
\qed
}\end{ex}

Integrals of this kind generalize the integral character formula due to
W.~Schmid and K.~Vilonen \cite{SchV2} for representations of $G_{\BB R}$
constructed from a $G_{\BB R}$-equivariant sheaf ${\cal F}$.
In that character formula the manifold $M$ is the flag variety ${\cal B}$ of
$\g g_{\BB C}$, $\Lambda = Ch({\cal F})$, and the integrand is the
pullback of a naturally defined form on a complex coadjoint orbit to
$T^*{\cal B}$ via the ``twisted moment map'' and can be be put into
the shape $\talpha$.

\separate

Let $T_{\BB C}$ be a maximal complex torus contained in $G_{\BB C}$,
i.e. $T_{\BB C}$ is a maximal subgroup of $G_{\BB C}$ isomorphic to
$\BB C^{\times} \times \dots \times \BB C^{\times}$.
Our last assumption on the action of $G_{\BB C}$ on $M$ is that
the points in $M$ fixed by the action of $T_{\BB C}$ are isolated.
Then there are only finitely many of them because $M$ is compact.
Since all the maximal tori of $G_{\BB C}$ are conjugate,
if this assumption holds for one torus $T_{\BB C}$ then it holds for all
maximal tori.

We denote by $\g g_{\BB C}^{rs}$ the set of {\em regular semisimple} elements
in $\g g_{\BB C}$.
These are elements $X \in \g g_{\BB C}$ such that the adjoint action of
$ad(X)$ on $\g g_{\BB C}$ is diagonalizable and has maximal possible rank.
We also denote by $\g g_{\BB R}^{rs} = \g g_{\BB R} \cap \g g_{\BB C}^{rs}$
the set of regular semisimple elements of $\g g_{\BB R}$.
It is an open and dense subset of $\g g_{\BB R}$.

For a regular semisimple element $X \in \g g_{\BB C}^{rs}$ we set
$\g t_{\BB C}(X) \subset \g g_{\BB C}$ to be the unique Cartan subalgebra
of $\g g_{\BB C}$ containing $X$ and $T_{\BB C}(X) = \exp (\g t_{\BB C}(X))$
to be the corresponding maximal torus.
Let $p \in M$ be a point fixed by $T_{\BB C}(X)$, then the complex Lie action
${\cal X} \mapsto {\cal L}(X_M){\cal X} =[X_M,{\cal X}]$
on the holomorphic vector fields ${\cal X}$ of $M$
gives rise to a linear transformation $L_p^{\BB C}$ on $T_pM$.
We define a function
$$
\den_p(X)= \det(L_p^{\BB C})
$$
which will appear in the denominator of the contribution of $p \in M_0(X)$
to the localization formula.

We will use the following description of $\den_p(X)$ which can serve
as an alternative definition.
The maximal torus $T_{\BB C}(X)$ acts linearly on $T_pM$.
Thus $T_pM$, as a representation of $T_{\BB C}(X)$,
decomposes into a direct sum of one-dimensional representations
$$
\BB C_{\beta_{p,1}} \oplus \dots \oplus \BB C_{\beta_{p,n}},
\qquad \beta_{p,1}, \dots, \beta_{p,n} \in \g t_{\BB C}(X)^*,
$$
where the action of $T_{\BB C}(X)$ on the one-dimensional complex vector
space $\BB C_{\beta_{p,k}}$ is given by
$$
\exp(Y) \cdot v = e^{\beta_{p,k}(Y)} v,
\qquad Y \in \g t_{\BB C}(X),\: v \in \BB C_{\beta_{p,k}}.
$$
The set of weights $\{ \beta_{p,1}, \dots, \beta_{p,n} \}$
is determined uniquely up to permutation.
Then we have
$$
\den_p(X)=\beta_{p,1}(X) \dots \beta_{p,n}(X).
$$
Notice that if the eigenvalues of $ad(X)$ are all purely imaginary
(that is $X$ lies in the Lie algebra of a compact subgroup of $G_{\BB C}$),
then we have the following relationship:
$$
\den_p(X) = i^n \cdot {\det}^{1/2}(L_p).
$$

%For a regular semisimple element $X \in \g g_{\BB C}^{rs}$ we set
%$\g t_{\BB C}(X) \subset \g g_{\BB C}$ to be the unique Cartan subalgebra
%of $\g g_{\BB C}$ containing $X$ and $T_{\BB C}(X) = \exp (\g t_{\BB C}(X))$
%to be the corresponding maximal torus.
%Let $p \in M$ be a point fixed by $T_{\BB C}(X)$, then
%$T_{\BB C}(X)$ acts linearly on $T_pM$.
%The representation of $T_{\BB C}(X)$ on $T_pM$ decomposes into a direct
%sum of one-dimensional representations
%$$
%\BB C_{\beta_{p,1}} \oplus \dots \oplus \BB C_{\beta_{p,n}},
%\qquad \beta_{p,1}, \dots, \beta_{p,n} \in \g t_{\BB C}(X)^*,
%$$
%where the action of $T_{\BB C}(X)$ on the one-dimensional complex vector
%space $\BB C_{\beta_{p,k}}$ is given by
%$$
%\exp(Y) \cdot v = e^{\beta_{p,k}(Y)} v,
%\qquad Y \in \g t_{\BB C}(X),\: v \in \BB C_{\beta_{p,k}}.
%$$
%The set of weights $\{ \beta_{p,1}, \dots, \beta_{p,n} \}$
%is determined uniquely up to permutation.
%We define a function
%$$
%\den_p(X)=\beta_{p,1}(X) \dots \beta_{p,n}(X)
%$$
%which will appear in the denominator of the contribution of $p \in M_0(X)$
%to the localization formula.
%
%Notice that if the eigenvalues of $ad(X)$ are all purely imaginary
%(that is $X$ lies in a Lie algebra of a compact subgroup of $G_{\BB C}$),
%then we have the following relationship:
%$$
%\den_p(X) = i^n \cdot {\det}^{1/2}(L_p).
%$$

We let $\Delta(X)$ denote the set of all weights that occur this way:
\begin{multline*}
\Delta(X) =
\{ \beta_{p,k} \in \g t_{\BB C}(X)^* ;\:
\text{$\beta_{p,k}$ appears in the weight decomposition}  \\
T_pM \simeq \BB C_{\beta_{p,1}} \oplus \dots \oplus \BB C_{\beta_{p,k}}
\oplus \dots \oplus \BB C_{\beta_{p,n}} \text{ for some $p \in M_0(X)$}\}.
\end{multline*}
It is a finite subset of $\g t_{\BB C}(X)^*  \setminus \{0\}$.

For instance, if $M$ is the flag variety of $\g g_{\BB C}$, then
$\Delta(X)$ is the root system of $\g g_{\BB C}(X)$ corresponding
to the choice of Cartan algebra $\g t_{\BB C}(X)$.

Let $\g g_{\BB R}'$ denote the set of regular semisimple elements
$X \in \g g_{\BB R}^{rs}$ which satisfy the following additional properties.
If $\g t_{\BB R}(X) \subset \g g_{\BB R}$ and
$\g t_{\BB C}(X) \subset \g g_{\BB C}$ are the unique Cartan
subalgebras in $\g g_{\BB R}$ and $\g g_{\BB C}$ respectively
containing $X$, then:
\begin{enumerate}
\item
The set of zeroes $M_0(X)$ is exactly the set of points in $M$ fixed by
the complex torus
$T_{\BB C}(X) = \exp (\g t_{\BB C}(X)) \subset G_{\BB C}$;
\item
$\beta(X) \ne 0$ for all $\beta \in \Delta(X) \subset \g t_{\BB C}(X)^*$;
\item
For each $\beta \in \Delta(X)$, we have either
\begin{equation}  \label{g'}
\re(\beta)|_{\g t_{\BB R}(X)} \equiv 0
\qquad \text{or} \qquad \re(\beta(X)) \ne 0.
\end{equation}
\end{enumerate}
Clearly, $\g g_{\BB R}'$ is an open subset of $\g g_{\BB R}$;
since $M$ is compact and $\Delta(X)$ is finite, the complement of
$\g g_{\BB R}'$ in $\g g_{\BB R}$ has measure zero;
and $\den_p(X) \ne 0$ for all $X \in \g g_{\BB R}'$.

\separate

The contribution to the integral of each zero $p \in M_0(X)$ will be
counted with some multiplicity $m_p \in \BB Z$ which we describe next.
We use the Bialynicki-Birula decomposition as restated in
Theorem 2.4.3 in \cite{ChG}.
Let $\BB C^{\times}$ be a subgroup of $G_{\BB C}$ such that the set of
fixed points $M^{\BB C^{\times}}$ in $M$ is finite. We embed $\BB C^{\times}$
into $\BB C$ in the most natural way so that
$\BB C^{\times} = \BB C \setminus \{0\}$.
For each fixed point $p \in M^{\BB C^{\times}}$
we define the {\em attracting} set
$$
O_p = \{ x \in M;\: \lim_{z \to 0} z^{-1} \cdot x =p \}.
$$
Clearly $p$ is the only point in $O_p$ fixed by $\BB C^{\times}$.
There is also a natural $\BB C^{\times}$-action on the tangent
space $T_pM$. It decomposes into a direct sum
\begin{equation}  \label{T-+}
T_pM = T_p^-M \oplus T_p^+M,
\end{equation}
$$
T_p^-M = \bigoplus_{k<0, \: k \in \BB Z} T_pM(k), \qquad
T_p^+M = \bigoplus_{k>0, \: k \in \BB Z} T_pM(k),
$$
where
$$
T_pM(k) =\{ v \in T_pM ;\:
z \cdot v = z^k v, \: \forall z \in \BB C^{\times} \}.
$$
Then we get the Bialynicki-Birula decomposition of $M$
into attracting sets $O_p$, each isomorphic to $T_p^-M$:

\begin{thm}(Bialynicki-Birula Decomposition \cite{BB})  \label{BB}

\noindent
\begin{enumerate}
\item
The attracting sets form a decomposition
$$
M = \coprod_{p \in M^{\BB C^{\times}}} O_p
$$
into smooth locally closed algebraic varieties;
\item
There are natural isomorphisms of algebraic varieties
\begin{equation}  \label{linear}
O_p \simeq T_p(O_p) \simeq T_p^-M
\end{equation}
which commute with the $\BB C^{\times}$-action.
\end{enumerate}
\end{thm}

\separate

Now let $X \in \g g_{\BB R}'$, and let $\g t_{\BB C}(X)$ and
$T_{\BB C}(X) = \exp(\g t_{\BB C}(X))$ be the corresponding complex
Cartan subalgebra and subgroup respectively.
Pick any $X' \in \g t_{\BB R}(X) \cap \g g_{\BB R}'$
in the same connected component of $\g t_{\BB R}(X) \cap \g g_{\BB R}'$
as $X$ and such that
$$
\re \beta(X) >0 \quad \Longleftrightarrow \quad \re \beta(X') >0
\qquad \text{and} \qquad
\re \beta(X) <0 \quad \Longleftrightarrow \quad \re \beta(X') <0
$$
for all $\beta \in \Delta(X)$, and the complex 1-dimensional subspace
$\{tX';\: t \in \BB C \} \subset \g g_{\BB C}$ is the Lie algebra of
a closed algebraic subgroup $\BB C^{\times}(X') \subset G_{\BB C}$
isomorphic to $\BB C^{\times}$.
Fix an isomorphism $\BB C^{\times}(X') \simeq \BB C^{\times}$ so that
the induced isomorphism of Lie algebras
$\{tX';\: t \in \BB C \} \simeq \BB C$ sends $X'$ into an element with
nonnegative real part.
We apply Theorem \ref{BB} to $\BB C^{\times}(X')$.
Then the set of points in $M$ fixed by $\BB C^{\times}(X')$ is just
$M_0(X') = M_0(X) = \{x_1,\dots,x_d\}$, say.
Let $O_k \subset M$ denote the attracting set of $x_k$ (instead of $O_{x_k}$).

For instance, if $M$ is the flag variety of $\g g_{\BB C}$, then the sets
$O_1,\dots,O_d$ are the orbits of a suitably chosen Borel subgroup containing
$T_{\BB C}(X)$, and the number of orbits $d$ equals the order of the
Weyl group of $\g g_{\BB C}$.

Since $\BB C^{\times}(X')$ is a subgroup of the torus $T_{\BB C}(X)$, their
actions commute, and the action of $T_{\BB C}(X)$ preserves each $O_k$.
Moreover, the proof of Theorem \ref{BB} shows that the isomorphism of
varieties (\ref{linear}) is $T_{\BB C}(X)$-equivariant.
In particular, the direct sum decomposition (\ref{T-+}) is a decomposition
of $T_{\BB C}(X)$-representations.

We define the multiplicity of a complex of sheaves ${\cal F} \in \C {M}$
at $x_k$ to be the Euler characteristic
\begin{equation}  \label{m}
m_k(X) = \chi \bigl( R\Gamma_{\{x_k\}}({\cal F}|_{O_k})_{x_k} \bigr)
= \chi \bigl( (j_{\{x_k\} \hookrightarrow O_k})^! ({\cal F}|_{O_k}) \bigr).
\end{equation}
The number $m_k(X)$ is an integer which
is exactly the local contribution of $x_k$ to the Lefschetz
fixed point formula, as generalized to sheaf cohomology by
M.~Goresky and R.~MacPherson \cite{GM}.

\separate

Now we are ready to state the main result of this article.

\begin{thm}  \label{main}
Let $G_{\BB C}$ act complex algebraically on a smooth complex projective
variety $M$ so that some (hence any) maximal torus
$T_{\BB C} \subset G_{\BB C}$ acts
with isolated fixed points. Suppose that a map
$\alpha : \g g_{\BB C} \to \Omega^*(M)$ satisfies Conditions \ref{conditions}.
And let $\Lambda \in {\cal L}^+_{G_{\BB R}} (M)$ be a $G_{\BB R}$-invariant
conic real-Lagrangian cycle in $T^*M$
such that the holomorphic moment map $\mu: T^*M \to \g g_{\BB C}^*$ is
proper on the set $\supp(\sigma|_{\Lambda})$.
Then, if $\phi$ is a smooth compactly supported differential
form on $\g g_{\BB R}'$ of top degree,
$$
\int_{\Lambda} \mu^* (\widehat{\phi\alpha}) \wedge e^{\sigma}
= \int_{\g g_{\BB R}} F_{\alpha} \phi,
$$
where $F_{\alpha}$ is an $Ad(G_{\BB R} \cap U_{\BB R})$-invariant
function on $\g g_{\BB R}'$ given by the formula
\begin{equation}  \label{mainequation}
F_{\alpha}(X) = (-2\pi i)^n
\sum_{k=1}^d m_k(X) \frac {\alpha(X)_{[0]}(x_k)}{\den_{x_k}(X)},
\end{equation}
where $n= \dim_{\BB C}(M)$, $\{x_1,\dots,x_d\} = M_0(X)$
is the set of zeroes of the vector field $X_M$ on $M$, and $m_k(X)$'s
are certain integer multiplicities.

To specify the multiplicities, let ${\cal F} \in \C (M)$
be a bounded complex of $G_{\BB R}$-equivariant sheaves on $M$ with
$\BB R$-constructible cohomology such that $Ch({\cal F}) = \Lambda$,
then the multiplicities are determined by the formula (\ref{m}).
%or a global formula
%\begin{equation}  \label{m2}
%m_k(X) = \chi (M, {\cal F}_{O_k}) =
%\chi \bigl( M, (j_{O_k \hookrightarrow M})_! \circ
%(j_{O_k \hookrightarrow M})^* ({\cal F}) \bigr).
%\end{equation}

We extend the function $F_{\alpha}$ by zero to a measurable function
on $\g g_{\BB R}$. If $F_{\alpha}$ happens to be locally integrable
with respect to the Lebesgue measure on
$\g g_{\BB R} \simeq \BB R^{\dim_{\BB R} \g g_{\BB R}}$,
then the equation (\ref{mainequation}) holds for smooth differential
forms $\phi$ of top degree which are compactly supported on $\g g_{\BB R}$
(and not necessarily on $\g g_{\BB R}'$).
\end{thm}

We divide the argument into two parts and give the proof in
sections \ref{deformationsection} and \ref{proof}.
We can say more about the multiplicities $m_k(X)$:

\begin{prop}  \label{m1=m2}
For each $X \in \g g_{\BB R}'$ and each bounded complex of
$G_{\BB R}$-equivariant sheaves ${\cal F} \in \C (M)$ with
$\BB R$-constructible cohomology, the multiplicities defined by the local
formula (\ref{m}) can also be given by a global formula
\begin{equation}  \label{m2}
m_k(X) = \chi (M, {\cal F}_{O_k}) =
\chi \bigl( M, (j_{O_k \hookrightarrow M})_! \circ
(j_{O_k \hookrightarrow M})^* ({\cal F}) \bigr).
\end{equation}
Moreover, these multiplicities depend on $Ch({\cal F})$ only and not on
the complex ${\cal F}$.
\end{prop}

\begin{rem}
In the special case when $\lambda$ equals $M$ as oriented cycles,
$Ch({\cal F})$ is $U_{\BB R}$-invariant, each multiplicity
$m_k(X)$ equals 1 and this theorem can be easily deduced from
the classical Berline-Vergne localization formula (Theorem \ref{BGV}).
\end{rem}

\begin{rem}
Notice that the cycle $\Lambda$ is invariant with respect to the
action of the group $G_{\BB R}$ which need not be compact, while the form 
$\alpha: \g g_{\BB C} \to \Omega^*(M)$ is required to be equivariant
with respect to a different group $U_{\BB R}$ only, and $U_{\BB R}$
may not preserve the cycle $\Lambda$.

The condition of the theorem that the moment map $\mu$ is proper on the set
$\supp(\sigma|_{\Lambda})$ is automatically satisfied when $\mu$
is proper on the support of the characteristic cycle $|\Lambda|$.
%It will be proved in Lemma \ref{m=m2} that the local formula (\ref{m})
%and the global formula (\ref{m2}) for the coefficient $m_k(X)$ agree
%with each other.
\end{rem}

\begin{rem}
Let ${\cal Z}({\cal U} (\g g_{\BB R}))$ denote the center of the
universal enveloping algebra of $\g g_{\BB R}$.
It is canonically isomorphic to the algebra of conjugate-invariant
constant coefficient differential operators on $\g g_{\BB R}$.
Suppose, in addition, that the distribution $\Delta$ on
$\g g_{\BB R}$ defined by
$$
\Delta: \phi \mapsto
\int_{\Lambda} \mu^* (\widehat{\phi\alpha}) \wedge e^{\sigma}
$$
is $Ad(G_{\BB R})$-invariant and is an eigendistribution for
${\cal Z}({\cal U} (\g g_{\BB R}))$ (i.e. each element of
${\cal Z}({\cal U} (\g g_{\BB R}))$ acts on $\Delta$ by
multiplication by some scalar).
Such situation arises in \cite{SchV2}, \cite{L1} and \cite{L2}
where the distribution $\Delta$ is the character of some virtual
representation of $G_{\BB R}$.
Then by Harish-Chandra's regularity theorem
(\cite{HC} or Theorem 3.3 in \cite{A}),
the function $F_{\alpha}$ from Theorem \ref{main}
is an $Ad(G_{\BB R})$-invariant, locally $L^1$ function
on $\g g_{\BB R}$ which is represented by a real analytic function on
the set of regular semisimple elements $\g g_{\BB R}^{rs}$.
Hence by the second part of Theorem \ref{main},
$$
\Delta(\phi) =
\int_{\Lambda} \mu^* (\widehat{\phi\alpha}) \wedge e^{\sigma}
= \int_{\g g_{\BB R}} F_{\alpha} \phi
$$
as distributions on $\g g_{\BB R}$.
\end{rem}

\separate

\begin{section}
{Deformation of $Ch({\cal F})$ in $T^*M$}  \label{deformationsection}
\end{section}

In this section ${\cal F} \in \C (M)$ is a bounded complex of
$G_{\BB R}$-equivariant sheaves on $M$ with $\BB R$-constructible cohomology
and $\Lambda = Ch({\cal F})$.
Recall that $B_R$ is an open ball in $\g g_{\BB C}^*$ defined by (\ref{B_R}).
We rewrite the integral (\ref{int}) as
\begin{multline}  \label{theintegral}
\int_{Ch({\cal F})}
\mu^* (\widehat {\phi\alpha}) \wedge e^{\sigma}
= \int_{Ch({\cal F})} \Bigl( \int_{\g g_{\BB R}} 
e^{\langle X, \mu(\zeta) \rangle} \phi(X) \wedge \alpha(X) \wedge e^{\sigma}
\Bigr) \\
=\lim_{R \to \infty}
\int_{\g g_{\BB R}' \times (Ch({\cal F}) \cap \mu^{-1}(B_R))}
e^{\langle X, \mu(\zeta) \rangle} \phi(X) \wedge \alpha(X) \wedge e^{\sigma}.
\end{multline}
(Of course, the orientation on
$\g g_{\BB R}' \times (Ch({\cal F}) \cap \mu^{-1}(B_R))$
is induced by the product orientation on
$\g g_{\BB R} \times Ch({\cal F})$.)
We will interchange the order of integration:
integrate over the characteristic cycle first and only then
perform integration over $\g g_{\BB R}'$.
By Lemma \ref{closed} the integrand in (\ref{theintegral})
is a closed differential form.

In this section we start with an element $X \in \g g_{\BB R}'$
and the characteristic cycle $Ch({\cal F})$
of a $G_{\BB R}$-equivariant complex of sheaves ${\cal F}$
on the projective variety $M$
and use general results of Section 4 in \cite{L1} to deform
$Ch({\cal F})$ into a cycle of the form
$$
m_1(X) T^*_{x_1}M +\dots+ m_d(X) T^*_{x_d}M,
$$
where $m_1(X),\dots,m_d(X)$ are the integer multiplicities given by
the equations (\ref{m}) and (\ref{m2}),
$x_1,\dots,x_d$ are the zeroes of the vector field $X_M$ on $M$,
and each cotangent space $T^*_{x_k}M$ is given some orientation.
Moreover, to ensure good behavior of our integral (\ref{theintegral}),
we will stay during the process of deformation inside the set
\begin{equation}  \label{deformationset}
\{\zeta \in T^*M ;\: \re( \langle X,\mu(\zeta) \rangle) \le 0\}.
\end{equation}
The precise result is stated in Proposition \ref{C(X)prop}.
This deformation will help us to calculate the integral (\ref{theintegral}).

\separate

Let $X \in \g g_{\BB R}'$, and let $\{x_1,\dots,x_d\} = M_0(X)$
be the set of zeroes of the vector field $X_M$ on $M$.
Let $\g t_{\BB C}(X) \subset \g g_{\BB C}$ and
$\g t_{\BB R}(X) \subset \g g_{\BB R}$ be the corresponding Cartan
subalgebras, and let
$T_{\BB C}(X) = \exp(\g t_{\BB C}(X)) \subset G_{\BB C}$ and
$T_{\BB R}(X) = \exp(\g t_{\BB R}(X)) \subset G_{\BB R}$ be the
corresponding connected subgroups.
Note that because we require $T_{\BB R}(X)$ to be connected it may not
be a Cartan subgroup of $G_{\BB R}$.

As a representation of $T_{\BB C}(X)$, the tangent space $T_{x_k}M$
at each zero $x_k$ decomposes into the direct sum (\ref{T-+}).
The space $T_{x_k}^-M$ in turn decomposes into a direct sum of
one-dimensional representations:
$$
T_{x_k}^-M \simeq 
\BB C_{\beta_{x_k,i_1}} \oplus \dots \oplus \BB C_{\beta_{x_k,i_m}}, \qquad
\{ \beta_{x_k,i_1}, \dots, \beta_{x_k,i_m} \} \subset
\{ \beta_{x_k,1}, \dots, \beta_{x_k,n} \}.
$$
By construction,
$$
\re \beta_{x_k,l}(X) <0 \quad \Rightarrow \quad
\beta_{x_k,l} \in \{ \beta_{x_k,i_1}, \dots, \beta_{x_k,i_m} \},
$$
$$
\re \beta_{x_k,l}(X) >0 \quad \Rightarrow \quad
\beta_{x_k,l} \not\in \{ \beta_{x_k,i_1}, \dots, \beta_{x_k,i_m} \}.
$$
Choose a linear coordinate $z_l: \BB C_{\beta_{x_k,l}} \, \tilde \to \, \BB C$
and define an inner product $\langle \cdot,\cdot\rangle_k$ on $T_{x_k}M$ by
$$
\bigl\langle(z_1,\dots,z_n), (z_1',\dots,z_n')\bigr\rangle_k =
z_1 \bar z_1' +\dots+ z_n \bar z_n'.
$$
Let $\|.\|_k$ be the respective norm on $T_{x_k}M$:
$$
\|(z_1,\dots,z_n)\|_k = |z_1|^2 +\dots+ |z_n|^2.
$$

Then, using the Bialynicki-Birula decomposition as stated in Theorem \ref{BB},
we obtain a decomposition of $M$ into smooth locally closed algebraic
varieties:
$$
M = \coprod_{k=1}^d O_k,
$$
where each $O_k$ is the attracting set of $x_k$, and we denote by
\begin{equation}  \label{psi}
\psi_{X,k}: T_{x_k}^-M \, \tilde \to \, O_k
\end{equation}
the $T_{\BB C}(X)$-equivariant isomorphism of varieties (\ref{linear}).

\separate

\begin{rem}  \label{triangle1}
Suppose ${\cal G}$ is a complex of sheaves on $M$ and $Z$ is a
locally closed subset of $M$. Let $i: Z \hookrightarrow M$ be
the inclusion. Then M. Kashiwara and P. Schapira introduce
in \cite{KaScha}, Chapter II, a complex $i_! \circ i^* ({\cal G})$
denoted by ${\cal G}_Z$.
If $Z'$ is closed in $Z$, then they prove existence of a
distinguished triangle
$$
{\cal G}_{Z \setminus Z'} \to {\cal G}_Z \to {\cal G}_{Z'}.
$$
Hence, by the additivity property of characteristic cycles,
$$
Ch( {\cal G}_Z ) = Ch( {\cal G}_{Z \setminus Z'} )
+ Ch( {\cal G}_{Z'}).
$$
\end{rem}

It follows that, as an element of $\K$ -- the Grothendieck group of
$\D(M)$, our complex of sheaves $\cal F$ is equivalent to
${\cal F}_{O_1}+ \dots + {\cal F}_{O_d}$, and so
$$
Ch({\cal F})=Ch({\cal F}_{O_1})+ \dots + Ch({\cal F}_{O_d}).
$$

The idea is to deform each summand $Ch({\cal F}_{O_k})$
separately.
Since $O_k$ is locally closed, there exists an open subvariety $U_k$ of $M$
containing $O_k$ as a closed subvariety.
Then by Proposition 4.22 of \cite{DM} or Section 4 of \cite{SchV1}
there exists a real-valued semi-algebraic ${\cal C}^2$-function
$f_k$ on $M$ such that $f_k$ is strictly positive on $U_k$ and
the boundary $\partial U_k$ is precisely the zero set of $f_k$.

\begin{lem}  \label{R}
There exists an $R>0$ such that, for each
$\zeta \in T_{x_k}^-M \subset T_{x_k}M$
with $\|\zeta\|_k \ge R$, the single-variable function
$$
f_k^{\zeta}(t) = f_k (\psi_{X,k}(t\zeta)), \qquad t \in \BB R,
$$
is strictly monotone decreasing for $t > 1/2$.
\end{lem}

\pf
Easily follows from the results on o-minimal structure described in \cite{DM},
and in particular the Monotonicity Theorem 4.1.
\qed

The dual space to $T_{x_k}^-M$, $(T_{x_k}^-M)^*$, can be regarded as a
subspace of the cotangent space at $x_k$:
$$
(T_{x_k}^-M)^* \quad \subset \quad T_{x_k}^*M
\quad = \quad (T_{x_k}^-M)^* \oplus (T_{x_k}^+M)^*.
$$
Let $B_k$ be the $\psi_{X,k}$-image of the open ball of radius $R$
$$
\{ \zeta \in (T_{x_k}^-M)^* ;\: \|\zeta\|_k < R \};
$$
$B_k$ is an open subset of $O_k$.

According to Remark \ref{triangle1} we have a distinguished triangle:
$$
{\cal F}_{B_k} \to {\cal F}_{O_k} \to {\cal F}_{O_k \setminus B_k},
$$
and hence
\begin{equation}  \label{tr1}
Ch({\cal F}_{O_k}) = Ch({\cal F}_{B_k}) + Ch({\cal F}_{O_k \setminus B_k}).
\end{equation}

Recall that the sheaf ${\cal F}$ is $G_{\BB R}$-equivariant.
In particular, $Ch({\cal F})$ is $T_{\BB R}(X)$-invariant, and so
$$
\re( \langle Y,\mu(\zeta) \rangle) =
- \re( \langle Y_M, \zeta \rangle) = 0
$$
for all $Y \in \g t_{\BB R}(X)$ and all $\zeta \in |Ch({\cal F})|$.

Similarly, because the set $O_k$ is $T_{\BB R}(X)$-invariant,
the sheaf ${\cal F}_{O_k}$ is $T_{\BB R}(X)$-equivariant too,
its characteristic cycle is $T_{\BB R}(X)$-invariant, and
$\re( \langle Y,\mu(\zeta) \rangle) =0$ for all 
$Y \in \g t_{\BB R}(X)$ and all $\zeta \in |Ch( {\cal F}_{O_k} )|$.

On the other hand, $B_k$ is an open subset of $O_k$ such that the vector
field $X_M$ is either tangent to the boundary $\partial B_k$ or
points outside $B_k$, but never points inside $B_k$.
It follows from the Open Embedding Theorem (Theorem \ref{oet})
and Lemma \ref{R} that
$\re( \langle Y,\mu(\zeta) \rangle) \le 0$ for all 
$Y \in \g t_{\BB R}(X)$ and all $\zeta \in |Ch( {\cal F}_{B_k} )|$.
Since
$Ch({\cal F}_{O_k \setminus B_k}) = Ch({\cal F}_{O_k}) - Ch({\cal F}_{B_k})$,
the same is true of $|Ch({\cal F}_{O_k \setminus B_k})|$.

\begin{lem}  \label{openembedding1}
The cycle $Ch({\cal F}_{O_k \setminus B_k})$
is homologous to the zero cycle inside the set
$$
\{\zeta \in T^*M ;\: \re (\langle X, \mu(\zeta) \rangle) \le 0 \}.
$$
\end{lem}

\pf
The sheaf ${\cal F}_{O_k \setminus B_k}$ is the extraordinary
direct image of a sheaf on $U_k$:
$$
{\cal F}_{O_k \setminus B_k} = (j_{U_k \hookrightarrow M})_! \circ
(j_{O_k \setminus B_k \hookrightarrow U_k})_! ({\cal F}|_{O_k \setminus B_k}).
$$

Recall that $f_k$ is real-valued semi-algebraic ${\cal C}^2$-function on $M$
which is strictly positive on $U_k$ and its zero set is precisely the
boundary $\partial U_k$.
It follows from the equation (\ref{mu}) and Lemma \ref{R} that,
for each $x \in O_k$ with $\| \psi^{-1}_{X,k}(x) \|_k > R/2$,
$$
\re\bigl(\langle X, \mu(df_k(x)) \rangle\bigr) =
-\re( \langle X_M, df_k(x) \rangle) \ge 0.
$$

By the Open Embedding Theorem (Theorem \ref{oet}),
\begin{multline*}
Ch( {\cal F}_{O_k \setminus B_k} ) =
Ch \bigl( (j_{U_k \hookrightarrow M})_! \circ 
(j_{O_k \setminus B_k \hookrightarrow U_k})_!
({\cal F}|_{O_k \setminus B_k}) \bigr)  \\
= \lim _{s \to 0^+}
Ch \bigl( (j_{O_k \setminus B_k \hookrightarrow U_k})_!
({\cal F}|_{O_k \setminus B_k}) \bigr) - s\frac {df_k}{f_k}.
\end{multline*}
Let $C$ be a $(2n+1)$-chain in $T^*M$
$$
C = - \Bigl ( Ch \bigl( (j_{O_k \setminus B_k \hookrightarrow U_k})_!
({\cal F}|_{O_k \setminus B_k}) \bigr) - s\frac {df_k}{f_k} \Bigr ),
\qquad s \in (0, \infty).
$$
Then $C$ is a conic chain,
its support $|C|$ lies inside the set (\ref{deformationset})
and the boundary of this chain $\partial C$ is 
$Ch( {\cal F}_{O_k \setminus B_k} )$ minus another cycle which we call
$$
\lim _{s \to +\infty}
Ch \bigl( (j_{O_k \setminus B_k \hookrightarrow U_k})_!
({\cal F}|_{O_k \setminus B_k}) \bigr) - s\frac {df_k}{f_k}.
$$
Notice that the last cycle is a cycle in $T^*M$ whose support lies
completely inside $T^*U_k$.

Recall the element $X' \in \g t_{\BB R} \cap \g g_{\BB R}'$
used to define attracting sets $O_1,\dots,O_d$.
Let $X'_{T_{x_k}M}$ be the vector field on $T_{x_k}M$ generated by $X'$.
Define a 1-form $\eta$ on $T_{x_k}^-M \setminus \{0\}$ to be
$$
\eta  =
\frac {\langle X'_{T_{x_k}M}, \: \cdot \: \rangle_k}
{\langle X'_{T_{x_k}M}, X'_{T_{x_k}M} \rangle_k}.
$$
We regard $\eta$ as a section of $T^*(O_k \setminus \{x_k\})$
via the isomorphism (\ref{psi}),
and let $\tilde \eta$ be any semi-algebraic extension of $\eta$
to a section of $T^*M|_{O_k \setminus \{x_k\}}$.
Since $X'$ lies in the same connected component of
$\g t_{\BB R}(X) \cap \g g_{\BB R}'$ as $X$,
it is easy to see that the real part
$\re \tilde \eta(X_M) = - \re \langle X, \mu(\tilde \eta) \rangle$
is strictly positive on $O_k \setminus \{x_k\}$.

Finally, define a $(2n+1)$-chain in $T^*M$
$$
\tilde C = - \Bigl(
\lim _{s \to +\infty}
Ch \bigl( (j_{O_k \setminus B_k \hookrightarrow U_k})_!
({\cal F}|_{O_k \setminus B_k}) \bigr) - s\frac {df_k}{f_k} \Bigr)
+ t \tilde \eta, \qquad t \in [0,\infty).
$$
Then its boundary
$$
\partial \tilde C =
\lim _{s \to +\infty}
Ch \bigl( (j_{O_k \setminus B_k \hookrightarrow U_k})_!
({\cal F}|_{O_k \setminus B_k}) \bigr) - s\frac {df_k}{f_k},
$$
$\tilde C$ is conic and its support $|\tilde C|$
lies in the set (\ref{deformationset}).
\qed

\separate

Next we deform $Ch({\cal F}_{B_k})$. We use another distinguished
triangle.

\begin{rem}
If ${\cal G}$ is a complex of sheaves on $M$,
$Z$ is a closed subset of $M$, $U=M \setminus Z$ is its complement
and $i: Z \hookrightarrow M$, $j:U \hookrightarrow M$ are the inclusion maps,
then we have a distinguished triangle
$$
(Ri)_* \circ i^! ({\cal G}) \to {\cal G} \to (Rj)_* \circ j^* ({\cal G}).
$$
Hence, by the additivity property of characteristic cycles,
$$
Ch( {\cal G} ) = Ch \bigl( (Ri)_* \circ i^! ({\cal G}) \bigr)
+ Ch \bigl( (Rj)_* \circ j^* ({\cal G}) \bigr).
$$
\end{rem}

We apply this remark with ${\cal G}={\cal F}_{B_k} =
(j_{B_k \hookrightarrow M})_! \circ (j_{B_k \hookrightarrow M})^* ({\cal F})$,
closed subset $Z = \{x_k\}$ and its complement $U=M \setminus \{x_k\}$:
\begin{multline}  \label{Ch(F_B)}
Ch({\cal F}_{B_k}) =  \\
Ch \bigl( (Rj_{\{x_k\} \hookrightarrow M})_* \circ
(j_{\{x_k\} \hookrightarrow M})^! ({\cal F}_{B_k}) \bigr)
+ Ch \bigl( (R j_{M \setminus \{x_k\} \hookrightarrow M})_* \circ
(j_{M \setminus \{x_k\} \hookrightarrow M})^* ({\cal F}_{B_k}) \bigr)
\end{multline}
Using that $B_k$ is an open subset of $O_k$, that the inclusion map
$O_k \hookrightarrow M$ is proper on the support of
$(j_{B_k \hookrightarrow O_k})_! ({\cal F}|_{B_k})$,
the Cartesian square
$$
\begin{matrix}
B_k & \hookrightarrow & O_k \\
\| & & \downarrow \\
B_k & \hookrightarrow & M
\end{matrix}
$$
and Proposition 3.1.9 of \cite{KaScha} we can write
\begin{multline*}
(j_{\{x_k\} \hookrightarrow M})^! ({\cal F}_{B_k})
= (j_{\{x_k\} \hookrightarrow M})^! \circ
(j_{B_k \hookrightarrow M})_! ({\cal F}|_{B_k})  \\
= (j_{\{x_k\} \hookrightarrow B_k})^! \circ
(j_{B_k \hookrightarrow M})^! \circ
(Rj_{O_k \hookrightarrow M})_* \circ
(j_{B_k \hookrightarrow O_k})_! ({\cal F}|_{B_k})  \\
= (j_{\{x_k\} \hookrightarrow B_k})^! \circ
(j_{B_k \hookrightarrow O_k})^* \circ
(j_{B_k \hookrightarrow O_k})_! ({\cal F}|_{B_k})
= (j_{\{x_k\} \hookrightarrow B_k})^! \circ ({\cal F}|_{B_k})  \\
= (j_{\{x_k\} \hookrightarrow B_k})^! \circ
(j_{B_k \hookrightarrow O_k})^! ({\cal F}|_{O_k})
= (j_{\{x_k\} \hookrightarrow O_k})^! ({\cal F}|_{O_k}).
\end{multline*}
Thus we can rewrite the equation (\ref{Ch(F_B)}) as
\begin{multline}  \label{tr2}
Ch({\cal F}_{B_k}) =  \\
Ch \bigl( (Rj_{\{x_k\} \hookrightarrow M})_* \circ
(j_{\{x_k\} \hookrightarrow O_k})^! ({\cal F}|_{O_k}) \bigr)
+ Ch \bigl( (R j_{M \setminus \{x_k\} \hookrightarrow M})_* \circ
(j_{M \setminus \{x_k\} \hookrightarrow M})^* ({\cal F}_{B_k}) \bigr).
\end{multline}
The cycle 
$Ch \bigl( (Rj_{\{x_k\} \hookrightarrow M})_* \circ
(j_{\{x_k\} \hookrightarrow O_k})^! ({\cal F}|_{O_k}) \bigr)$
is the cotangent space $T^*_{x_k} M$ equipped with
orientation (\ref{orientation}) and multiplicity $m_k(X)$ given by
the local formula (\ref{m}).

\separate

It remains to show that the second summand of (\ref{tr2}) is homologous
to zero.
Let ${\cal G}$ denote the sheaf
$(j_{M \setminus \{x_k\} \hookrightarrow M})^* ({\cal F}_{B_k})$
on $M \setminus \{x_k\}$; it is supported inside the closure of
$B_k \setminus \{x_k\}$ in $M \setminus \{x_k\}$.
Pick any real-valued semi-algebraic ${\cal C}^2$-function $\tilde f_k$ on $M$
such that $\tilde f_k$ is strictly positive on $M \setminus \{x_k\}$ and
$\tilde f(x_k)=0$. Similarly to Lemma \ref{R} we have:

\begin{lem}  \label{R'}
There exists an $R'>0$ such that, for each
$\zeta \in T_{x_k}^-M \subset T_{x_k}M$
with $\|\zeta\|_k \le R'$, the single-variable function
$$
\tilde f_k^{\zeta}(t) = \tilde f_k (\psi_{X,k}(t\zeta)), \qquad t \in \BB R,
$$
is strictly monotone increasing for $t \in [0,2]$.
\end{lem}

Since we are free to modify $\tilde f_k$ on any compact subset of $M$ which
does not contain $x_k$, we may assume that $R' > R$.

\begin{lem}  \label{openembedding2}
The cycle
$$
Ch \bigl( (R j_{M \setminus \{x_k\} \hookrightarrow M})_* \circ
(j_{M \setminus \{x_k\} \hookrightarrow M})^* ({\cal F}_{B_k}) \bigr)
= Ch \bigl( (R j_{M \setminus \{x_k\} \hookrightarrow M})_* ({\cal G}) \bigr)
$$
is homologous to the zero cycle inside the set
$$
\{\zeta \in T^*M ;\: \re (\langle X, \mu(\zeta) \rangle) \le 0 \}.
$$
\end{lem}

\pf
Except for a few obvious modifications, this proof is identical
to the proof of Lemma \ref{openembedding1}.
First we observe that because $R'$ from Lemma \ref{R'} is bigger
than $R$ used to define the set $B_k$,
for each $x \in \overline{B_k} \setminus \{x_k\}$,
we have $\re\bigl(\langle X, \mu(d \tilde f_k(x)) \rangle\bigr) =
-\re( \langle X_M, d \tilde f_k(x) \rangle) \le 0$.

By the Open Embedding Theorem (Theorem \ref{oet}),
$$
Ch \bigl( (R j_{M \setminus \{x_k\} \hookrightarrow M})_* ({\cal G}) \bigr)  \\
= \lim _{s \to 0^+}
Ch ({\cal G}) + s\frac {d \tilde f_k}{ \tilde f_k}.
$$
Thus we introduce a $(2n+1)$-chain in $T^*M$
$$
C' = - \Bigl ( Ch ({\cal G}) + s\frac {d \tilde f_k}{ \tilde f_k} \Bigr),
\qquad s \in (0, \infty).
$$
Then $C'$ is a conic chain,
its support $|C'|$ lies inside the set (\ref{deformationset})
and the boundary of this chain $\partial C'$ is 
$Ch \bigl( (R j_{M \setminus \{x_k\} \hookrightarrow M})_* ({\cal G}) \bigr)$
minus another cycle which we call
$$
\lim _{s \to +\infty}
Ch ({\cal G}) + s\frac {d \tilde f_k}{ \tilde f_k}.
$$
Notice that the last cycle is a cycle in $T^*M$ whose support lies
completely inside $T^*(M \setminus \{x_k\})$.

Recall the section $\tilde \eta$ of $T^*M|_{O_k \setminus \{x_k\}}$
constructed in the proof of Lemma \ref{openembedding1}.
It has the property that $\re \langle X, \mu(\tilde \eta) \rangle$ 
is strictly negative on $O_k \setminus \{x_k\}$.

Finally, define a $(2n+1)$-chain in $T^*M$
$$
\tilde C' = - \Bigl( \lim _{s \to +\infty}
Ch ({\cal G}) + s\frac {d \tilde f_k}{ \tilde f_k} \Bigr)
+ t \tilde \eta, \qquad t \in [0,\infty).
$$
Then its boundary
$$
\partial \tilde C' =
\lim _{s \to +\infty}
Ch ({\cal G}) + s\frac {d \tilde f_k}{ \tilde f_k},
$$
$\tilde C'$ is conic and its support $|\tilde C'|$ lies in
the set (\ref{deformationset}).
\qed

\separate

Combining the equations (\ref{tr1}), (\ref{tr2}) and
lemmas \ref{openembedding1}, \ref{openembedding2} we obtain
the following key result.

\begin{prop}  \label{C(X)prop}
For each element $X \in \g g_{\BB R}'$, there is a
Borel-Moore chain $C(X)$ in $T^*M$
of dimension $(2n+1)$ with the following properties:
\begin{enumerate}
\item
$C(X)$ is conic, i.e. invariant under the scaling action of the
multiplicative group of positive reals $\BB R^{>0}$ on $T^*M$;
\item
The support of $C(X)$ lies in the set
$\{\zeta \in T^*M;\: \re( \langle X, \mu(\zeta) \rangle ) \le 0 \}$;
\item
Let $x_1, \dots, x_d$ be the zeroes of the vector field $X_M$ on $M$, then
$$
\partial C(X) =
Ch({\cal F}) - \bigl( m_1(X) T^*_{x_1}M + \dots + m_d(X) T^*_{x_d}M \bigr),
$$
where $m_1(X),\dots,m_d(X)$ are the integer multiplicities determined
by the local formula (\ref{m}) 
and the orientation of $T^*_{x_k}M$ is chosen so that
if we write each $z_l$ as $x_l + iy_l$, then the
$\BB R$-basis\footnote{The holomorphic cotangent bundle $T^*M$ and
the ${\cal C}^{\infty}$ cotangent bundle $T^*M^{\BB R}$ are
identified according to Remark \ref{T*M}.}
\begin{equation}   \label{orientation}
\{ dx_1, dy_1, \dots, dx_n, dy_n \}  \qquad \text{of} \qquad
T^*_{x_k}M \simeq
(\BB C_{\beta_{x_k,1}})^* \oplus \dots \oplus (\BB C_{\beta_{x_k,n}})^*
\end{equation}
is positively oriented;
\item
Moreover, if $\tilde X \in \g t_{\BB R}(X) \cap \g g_{\BB R}'$
lies in the same connected component of $\g t_{\BB R}(X) \cap \g g_{\BB R}'$
as $X$, then the same choice of element
$X' \in \g t_{\BB R}(\tilde X) \cap \g g_{\BB R}'$ works for $\tilde X$.
In this case the chain $C(\tilde X)$ is identical to $C(X)$.
\end{enumerate}
\end{prop}

\begin{rem}  \label{orientationremark}
The holomorphic cotangent space $T^*_{x_k}M$ has a natural orientation
coming from its complex structure. This orientation need not agree with
the orientation given by (\ref{orientation}). In fact,
$$
\text{the complex orientation of $T^*_{x_k}M$}
= (-1)^n \text{ the orientation given by (\ref{orientation})}.
$$
\end{rem}

Next we show that the local formula (\ref{m}) and the global
formula (\ref{m2}) for the coefficient $m_k(X)$ give the same answer.

%\begin{lem}  \label{m=m2}
%For each $X \in \g g_{\BB R}'$, the right hand sides of the equations
%(\ref{m}) and (\ref{m2}) are equal.
%\end{lem}

\noindent {\it Proof of Proposition \ref{m1=m2}.}
By a generalization of the Hopf Index Theorem (Theorem \ref{hit}),
$$
\chi(M, {\cal F}_{O_k}) = \# \bigl( [M] \cap Ch({\cal F}_{O_k}) \bigr).
$$
%(Alternatively one can apply the equation (5.30) from \cite{Schu}.)
Then by Proposition \ref{C(X)prop}, the characteristic cycle
$Ch({\cal F}_{O_k})$ is homologous to the cycle
$$
\chi \bigl( R\Gamma_{\{x_k\}}({\cal F}|_{O_k})_{x_k} \bigr)
\cdot T^*_{x_k}M,
$$
where $T^*_{x_k}M$ is given orientation as described in (\ref{orientation}).
Since $T^*_{x_k}M$ intersects $M$ transversally, we see that the
right hand side of (\ref{m2}) is
$$
\chi(M, {\cal F}_{O_k}) = 
\# \bigl( [M] \cap Ch({\cal F}_{O_k}) \bigr) =
\chi \bigl( R\Gamma_{\{x_k\}}({\cal F}|_{O_k})_{x_k} \bigr).
\qquad \square
$$

\separate

\begin{section}
{Proof of Theorem \ref{main}}  \label{proof}
\end{section}

In this section we compute the integral (\ref{theintegral})
first under the assumption that the form $\phi$ is compactly
supported in $\g g_{\BB R}'$ and then in general.
First we define a deformation $\Theta_t(X): T^*M \to T^*M$,
where $X \in \g g_{\BB R}'$, $t \in [0,1]$.
It has the following purpose. In the classical proof of
the Fourier inversion formula
$$
\phi(X)= \frac 1{(2\pi i)^{\dim_{\BB R}\g g_{\BB R}}}
\int_{\xi \in i \g g_{\BB R}^*}
\hat \phi(\xi) e^{-\langle X,\xi \rangle}
$$
we multiply the integrand by a term like $e^{-t \|\xi\|^2}$
to make it integrable over $\g g_{\BB R} \times i \g g_{\BB R}^*$,
and then let $t \to 0^+$.
The deformation $\Theta_t(X)$ has a very similar effect --
it makes our integrand an $L^1$-object.
Proposition \ref{slanting} says that this substitution is permissible.
Its proof is very technical, but the idea is quite simple.
The difference between the original integral (\ref{theintegral})
and the deformed one is expressed by an integral of 
$e^{\langle X, \mu(\zeta) \rangle} \phi(X) \wedge \alpha(X) \wedge e^{\sigma}$
over a certain cycle $\tilde C(R)$ supported in
$\g g_{\BB R}' \times (T^*M \cap \{\|\mu(\zeta)\|_{\g g_{\BB C}^*} = R\})$
which depends on $R$ by scaling along the fiber.
Recall that the Fourier transform $\hat \phi$ decays rapidly
in the imaginary directions which is shown by integration
by parts. We modify this integration by parts argument to prove
a similar statement about the behavior of the integrand on
the support of $\tilde C(R)$ as $R \to \infty$.
Hence the difference of integrals in question tends to zero.

\separate

Pick an element $X_0$ lying in the support of $\phi$ and let
$\g t_{\BB R}(X_0) \subset \g g_{\BB R}$ be the Cartan subalgebra
containing $X_0$.
There exists an open neighborhood $\Omega$ of $X_0$ in
$\g g_{\BB R}'$ and a smooth map $\omega: \Omega \to G_{\BB R}$
with the following three properties:
\begin{enumerate}
\item
$\omega|_{\Omega \cap \g t_{\BB R}(X_0)} \equiv e$, the identity element of
$G_{\BB R}$;
\item
For every $X \in \Omega$,
the conjugate Cartan subalgebra
$\omega(X) \g t_{\BB R}(X_0) \omega(X)^{-1}$
contains $X$;
\item
$\omega(X) = \omega(Y)$ whenever $X, Y \in \Omega$ and
$\g t_{\BB C}(X) = \g t_{\BB C}(Y)$ (i.e. $[X,Y]=0$).
\end{enumerate}
Notice that if $X \in \Omega$, then $\omega(X) \cdot M_0(X_0) = M_0(X)$.
Making $\Omega$ smaller if necessary, we can assume that both
$\Omega$ and $\Omega \cap \g t_{\BB R}(X_0)$ are connected.
Let $\g t_{\BB C}(X_0) = \g t_{\BB R}(X_0) \oplus i \g t_{\BB R}(X_0)
\subset \g g_{\BB C}$
be the complex Cartan subalgebra containing $X_0$.

\begin{rem}
One cannot deal with the integral (\ref{theintegral})
``one Cartan algebra at a time'' and avoid introducing
a map like $\omega$ because the limit
$$
\lim_{R \to \infty}
\int_{\g t_{\BB R} \times (Ch({\cal F}) \cap \mu^{-1}(B_R))}
e^{\langle X, \mu(\zeta) \rangle} \phi(X) \wedge \alpha(X) \wedge e^{\sigma}.
$$
may not exist.
(Recall that $B_R$ is an open ball in $\g g_{\BB C}^*$ defined by (\ref{B_R}).)
\end{rem}

\separate

From now on we assume that the support of $\phi$ lies in
$\Omega$. The general case when $\supp(\phi) \subset \g g_{\BB R}'$
can be reduced to this special case
by a partition of unity argument.

\separate

Our biggest obstacle to making any deformation argument computing the
integral (\ref{theintegral}) is that the integration takes place over
a cycle which is not compactly supported and Stokes' theorem no longer applies.
In order to overcome this obstacle, we construct a deformation
$\Theta_t: \Omega \times T^*M \to \Omega \times T^*M$, $t \in [0,1]$,
such that $\Theta_0$ is the identity map;
$$
\re \bigl( (\Theta_t)^*\langle X, \mu(\zeta) \rangle \bigr)
<  \re(\langle X, \mu(\zeta) \rangle)
$$
for $t>0$, $X \in \Omega$ and $\zeta \in T^*M$
which does not lie in the zero section (Lemma \ref{realpartlemma});
$\Theta_t$ essentially commutes with scaling the fiber of $T^*M$
(Lemma \ref{scaling}).
The last two properties will imply that the integral
$$
\int_{\g g_{\BB R} \times (Ch({\cal F}) \cap \mu^{-1}(B_R))}
(\Theta_t)^* \bigl(
e^{\langle X, \mu(\zeta) \rangle} \phi(X) \wedge \alpha(X) \wedge e^{\sigma}
\bigr)
$$
converges absolutely for $t \in (0,1]$.
Finally, the most important property of $\Theta_t$ is stated
in Proposition \ref{slanting} which essentially says that we can
replace our integrand
$$
e^{\langle X, \mu(\zeta) \rangle} \phi(X) \wedge \alpha(X) \wedge e^{\sigma}
$$
with the pullback
$$
(\Theta_t)^* \bigl(
e^{\langle X, \mu(\zeta) \rangle} \phi(X) \wedge \alpha(X) \wedge e^{\sigma}
\bigr).
$$

\separate

We restate Theorem 1 of \cite{Su}:
\begin{prop}  \label{embedding}
There is a projective embedding $\nu: M \to \BB CP^N$ and
a group representation $\rho: G_{\BB C} \to PGL(N)$ such that
$\rho(g) \cdot \nu(x) = \nu(g \cdot x)$
for every $g \in G_{\BB C}$ and $x \in M$.
\end{prop}

Let $\{x_1, \dots, x_d \}$ be the set of zeroes $M_0(X_0)$.
For $D>0$, we denote by $B_D$ the open ball in $\BB C^n$ of radius $D$:
$$
B_D = \{(z_1,\dots,z_n) \in \BB C^n ;\: |z_1|^2+\dots+|z_n|^2 < D^2 \}.
$$
Using Proposition \ref{embedding} one can construct a ${\cal C}^{\infty}$
diffeomorphism onto an open subset $V_k \subset M$ containing $x_k$
$$
\psi_{X_0,k}: B_{4D} \, \tilde \to \, V_k
$$
such that $\psi_{X_0,k}(0) = x_k$ and, for each $Y \in \g t_{\BB C}(X_0)$,
the tangent map $d \psi_{X_0,k}$ sends the vector field on
$\BB C^n \simeq T_0B_{4D}$
$$
\beta_{x_k,1}(Y) z_1 \frac{\partial}{\partial z_1} + \dots +
\beta_{x_k,n}(Y) z_n \frac{\partial}{\partial z_n}
$$
into $-Y_M$.
Note that this condition implies $V_k \cap M_0(X_0) = \{x_k\}$.

On the other hand, each point $x \in M \setminus M_0(X_0)$ has a
${\cal C}^{\infty}$ chart
$$
\psi_{X_0,x}: B_{4D} \, \tilde \to \, V_x
$$
such that $\psi_{X_0,x}(0) = x$ and
\begin{equation} \label{1}
d \psi_{X_0,x} \Bigl( \frac{\partial}{\partial z_1} \Bigr) = - (X_0)_M.
\end{equation}
Making $V_x$ smaller if necessary, we can assume that
$V_x \cap M_0(X_0) = \varnothing$.
For $Y \in \Omega$, let
$$
Y_1^x(z) \frac{\partial}{\partial z_1} + \dots +
Y_n^x(z) \frac{\partial}{\partial z_n}
$$
be the inverse image of the vector field $-Y_M$
under the tangent map $d \psi_{X_0,x}$.
By continuity (\ref{1}) implies that there is an open neighborhood
$\Omega_x$ of $X_0$ such that $\re (Y_1^x(z)) > 0$ for $z \in B_{4D}$ and
$Y \in \Omega_x \cap \g t_{\BB C}(X_0)$.

We extend $\{ \psi_{X_0,1},\dots,\psi_{X_0,d} \}$ to an atlas
$\{ \psi_{X_0,1},\dots,\psi_{X_0,d'} \}$ of $M$ so that, for $d < k \le d'$,
$\psi_{X_0,k} = \psi_{X_0,x_k'}$ for some $x_k' \in M \setminus M_0(X_0)$
and the smaller open sets completely cover $M$:
\begin{equation}  \label{D}
\bigcup_{k=1}^{d'} \psi_{X_0,k} (B_D)=M.
\end{equation}
Set $V_k = \psi_{X_0,k} (B_{4D})$, $k=1,\dots,d'$.

For each $X \in \Omega$ we define maps
$$
\psi_{X,k}: B_{4D} \, \tilde \to \, \omega(X) \cdot V_k,
\qquad \psi_{X,k} (z) = \omega(X) \cdot \psi_{X_0,k}(z),
\qquad 1 \le k \le d'.
$$
Then $\{ \psi_{X,1},\dots,\psi_{X,d'} \}$ form another atlas of $M$.
Note that, for $k=1,\dots,d$, $\psi_{X,k}(0) = \omega(X) \cdot x_k$ and,
for each $Y \in \g t_{\BB C}(X) = \omega(X) \g t_{\BB C}(X_0) \omega(X)^{-1}$,
the tangent map $d \psi_{X,k}$ sends the vector field on $T_{x_k}M$
$$
\beta_{x_k,1} \bigl( \omega(X)^{-1} Y \omega(X) \bigr)
z_1 \frac{\partial}{\partial z_1} + \dots +
\beta_{x_k,n} \bigl( \omega^{-1}(X) Y \omega(X) \bigr)
z_n \frac{\partial}{\partial z_n}
$$
into $-Y_M$. We extend
$\beta_{x_k,1}, \dots, \beta_{x_k,n} \in \g t_{\BB C}(X_0)$ to $\Omega$ by
$$
\beta_{x_k,l} (Y) \, =_{def} \,
\beta_{x_k,l} \bigl( \omega(Y)^{-1} Y \omega(Y) \bigr),
\qquad Y \in \Omega, \: l=1,\dots,n.
$$
This way, for all $X \in \Omega$ and all $Y \in \g t_{\BB C}(X)$,
we can write
\begin{equation}  \label{vfield}
d \psi_{X,k} \Bigl(
\beta_{x_k,1} (Y) z_1 \frac{\partial}{\partial z_1} + \dots +
\beta_{x_k,n} (Y) z_n \frac{\partial}{\partial z_n} \Bigr)
= -Y_M.
\end{equation}
If $k=d+1,\dots,d'$ and $Y \in \Omega$, let
$$
Y_1^k(z) \frac{\partial}{\partial z_1} + \dots +
Y_n^k(z) \frac{\partial}{\partial z_n}
$$
be the inverse image of the vector field $-Y_M$
under the tangent map $d \psi_{Y,k}$.
Note that
$$
d \psi_{X,k} \Bigl( \frac{\partial}{\partial z_1} \Bigr) = - X_M,
\qquad d < k \le d'.
$$
Hence making $\Omega$ smaller if necessary, we can assume that
$\re (Y_1^k(z)) > 0$ for $d < k \le d'$, $z \in B_{4D}$,
$Y \in \Omega \cap \g t_{\BB C}(X)$ and all $X \in \Omega$.

Finally, we define maps
$$
\psi_k: \Omega \times B_{4D} \to \Omega \times M,
$$
$$
\psi_k(X,z) = (X, \psi_{X,k}(z)) = (X, \omega(X) \cdot \psi_{X_0,k}(z)),
\qquad 1 \le k \le d'.
$$
Each $\psi_k$ is a diffeomorphism onto its image, and their images for
$k=1,\dots,d'$ cover all of $\Omega \times M$.
Thus we obtain an atlas $\{ \psi_1,\dots,\psi_{d'} \}$ of $\Omega \times M$.

Expand $(z_1, \dots , z_n)$ to a standard coordinate system
$(z_1, \dots, z_n,\xi_1, \dots , \xi_n)$ on the cotangent space $T^* B_{4D}$
so that every element of $T^* B_{4D} \simeq B_{4D} \times \BB C^n$
is expressed in these coordinates as
$$
(z_1,\dots, z_n, \xi_1 dz_1 + \dots + \xi_n dz_n).
$$
This gives us a chart
$$
\tilde \psi_k: (X, z_1,\dots,z_n,\xi_1,\dots,\xi_n)
\to \Omega \times T^*M
$$
and an atlas $\{ \tilde \psi_1,\dots,\tilde \psi_{d'} \}$
of $\Omega \times T^*M$.
For $(z_1,\dots,z_n,\xi_1,\dots,\xi_n) \in T^*B_{4D}$,
define norms $\|z\| = \sqrt{|z_1|^2+\dots+|z_n|^2}$ and
$\|\xi\| = \sqrt{|\xi_1|^2+\dots+|\xi_n|^2}$.

Find an $\epsilon>0$ small enough so that for each $k=1,\dots,d$
\begin{equation} \label{one}
\psi_{X_0,k} (B_{\epsilon}) \cap
\overline{\bigcup_{l \ne k} \psi_{X_0,l} (B_{3D})} = \varnothing;
\end{equation}
we also assume that $\epsilon \le D/2$.

\separate

Since $G_{\BB R}$ acts on $M$ by complex automorphisms,
the symplectic form $\sigma$ in these coordinates is
$d\xi_1 \wedge dz_1 + \dots +d\xi_n \wedge dz_n$.

For $k=1,\dots,d$, the equations (\ref{mu}) and (\ref{vfield}) say that
the exponential part in the chart $\tilde \psi_k$ becomes
\begin{equation}  \label{exponent}
\langle X, \mu(\zeta) \rangle =
\beta_{x_k,1}(X)z_1\xi_1 + \dots + \beta_{x_k,n}(X)z_n\xi_n.
\end{equation}

Let $\delta: \BB R \to [0,1]$ be a smooth bump function which
takes on value $1$ on $[-D, D]$, vanishes outside
$(-2D,2D)$, and is nondecreasing on negative reals,
non-increasing on positive reals.
By making $\epsilon$ smaller if necessary we may assume that
$|2 \epsilon \delta'(x)|<1$ for all $x$.

Let $\gamma: \BB R^+ \to (0,1]$ be another smooth function
which is non-increasing, $\gamma([0,1])=\{1\}$,
$\gamma(x)=\frac 1x$ for $x >2$, and
$\frac 1x \le \gamma(x) \le \frac 2x$ for all $x \ge 1$.

And let $\rho: \BB R^+ \to [0,\infty)$ be a smooth monotone increasing function
such that its derivative $\rho'(x) \le \frac 12$ for all $x$ and
$$
\begin{cases}
\rho(x) = \frac 14 x^2  & \text{if $x \in [0,1]$;}  \\
\rho(x) = ax & \text{if $x \ge 2$}
\end{cases}
$$
for some constant $a>0$.

Note that the derivatives of $\delta$, $\gamma$ and $\rho$ are uniformly
bounded on their respective domains.

\separate

For each $t \in [0,1]$ and $k=1,\dots,d'$ we define a map
$\Theta_t^k: \Omega \times T^*M \to \Omega \times T^*M$.
If $Y \in \g t_{\BB C}(X_0)$ and $k =1,\dots,d$,
we define a diffeomorphism $\tilde \Theta_{Y,t}^k$
on $T^*B_{4D} \simeq B_{4D} \times \BB C^n$ by
\begin{align*}
\tilde \Theta_{Y,t}^k (z_1, \dots, z_n,\xi_1, \dots, \xi_n) &= 
(z_1', \dots, z_n',\xi_1, \dots, \xi_n)  \\
z_j' &= z_j - \frac {\overline{\beta_{x_k,j}(Y)}}{|\beta_{x_k,j}(Y)|}
t \epsilon \delta(\|z\|)\gamma(t \|\xi\|) \bar \xi_j,
\qquad 1 \le j \le d,
\end{align*}
(the requirement $|2 \epsilon \delta'|<1$ ensures that
$\tilde \Theta_{Y,t}^k$ is one-to-one).
%\begin{multline*}
%(z_1, \dots, z_n,\xi_1, \dots, \xi_n) \mapsto \\
%\Bigl( z_1 - \frac {\overline{\beta_{x_k,1}(Y)}}{|\beta_{x_k,1}(Y)|}
%t \epsilon \delta(\|z\|)\gamma(t \|\xi\|) \bar \xi_1, \dots,
%z_n - \frac{\overline{\beta_{x_k,n}(Y)}}{|\beta_{x_k,n}(Y)|}
%t \epsilon \delta(\|z\|)\gamma(t \|\xi\|) \bar \xi_n, \\
%\xi_1, \dots, \xi_n \Bigr).
%\end{multline*}
If $k=d+1,\dots,d'$, we define a diffeomorphism
$\tilde \Theta_{Y,t}^k$ on $T^*B_{4D} \simeq B_{4D} \times \BB C^n$ by
$$
\tilde \Theta_{Y,t}^k (z_1, \dots, z_n,\xi_1, \dots, \xi_n) = 
(z_1,\dots, z_n,
\xi_1 - \delta(\|z\|) \rho(t\|\xi\|), \xi_2,\dots, \xi_n)
$$
(again, the requirement $\rho' \le \frac 12$ ensures that
$\tilde \Theta_{Y,t}^k$ is one-to-one).
The map $\tilde \Theta_{Y,t}^k$ shifts $(z_1,\dots,z_n)$ by a vector
$$
- t \epsilon \delta(\|z\|)\gamma(t \|\xi\|)
\Bigl( \frac {\overline{\beta_{x_k,1}(Y)}}{|\beta_{x_k,1}(Y)|} \bar \xi_1,
\dots,
\frac {\overline{\beta_{x_k,n}(Y)}}{|\beta_{x_k,n}(Y)|} \bar \xi_n \Bigr)
\qquad \text{if $k=1,\dots,d$}
$$
which has length at most $2\epsilon \le D$ (because $\gamma(x) \le \frac 2x$),
and shifts $\xi_1$ by a scalar
$$
- \delta(\|z\|) \rho(t\|\xi\|) \qquad \text{if $k=d+1,\dots,d'$}.
$$
Hence the maps $\tilde \Theta^k_{Y,t}$ and $(\tilde \Theta_{Y,t}^k)^{-1}$
leave points outside the set $\{(z,\xi);\: \|z\| \le 2D\}$
completely unaffected.
Then we use the diffeomorphism between
$B_{4D} \times \BB C^n$ and $T^*V_k \subset T^*M$ induced by the map
$\psi_{X_0,k}: B_{4D} \to M$
to regard $\tilde \Theta_{Y,t}^k$ as a map on $T^*V_k$.
But since $\tilde \Theta_{Y,t}^k$ becomes the identity map when the
basepoint of $\zeta \in T^*M$ lies away from the compact subset
$$
\psi_{X_0,k} (\{z;\: \|z\| \le 2D \}) \subset V_k \subset M,
$$
$\tilde \Theta_{Y,t}^k$
can be extended by identity to a diffeomorphism $T^*M \to T^*M$.

Finally, we define
$\Theta_t^k: \Omega \times T^*M \to \Omega \times T^*M$
using the ``twisted'' product structure of $\Omega \times T^*M$
induced by $\omega(X)$.
Recall that the group $G_{\BB C}$ acts on $M$ which induces
an action on $T^*M$.
For $g \in G_{\BB C}$ and $\zeta \in T^*M$, we denote this action by
$g \cdot \zeta$.
Then, for $(X,\zeta) \in \Omega \times T^*M$, we set
$$
\Theta_t^k (X,\zeta) = \bigl(
X, \omega(X) \cdot (\tilde \Theta_{Y,t}^k (\omega(X)^{-1} \cdot \zeta)) \bigr),
\qquad \text{where $Y = \omega(X)^{-1} X \omega(X) \in \g t_{\BB C}(X_0)$.}
$$

Inside the chart $\tilde \psi_k$ centered at the point $(X_0,x_k)$,
$\Theta_t^k$ is formally given by the same expression as before:
\begin{align*}
\Theta_t^k (X, z_1, \dots, z_n,\xi_1, \dots, \xi_n) &=
(X, z_1', \dots, z_n',\xi_1, \dots, \xi_n)  \\
z_j' &= z_j - \frac {\overline{\beta_{x_k,j}(X)}}{|\beta_{x_k,j}(X)|}
t \epsilon \delta(\|z\|)\gamma(t \|\xi\|) \bar \xi_j,
\qquad 1 \le j \le d,
\end{align*}
%\begin{multline*}
%(X, z_1, \dots, z_n,\xi_1, \dots, \xi_n) \mapsto \\
%\Bigl( X, z_1 - \frac {\overline{\beta_{x_k,1}(X)}}{|\beta_{x_k,1}(X)|}
%t \epsilon \delta(\|z\|)\gamma(t \|\xi\|) \bar \xi_1, \dots,
%z_n - \frac {\overline{\beta_{x_k,n}(X)}}{|\beta_{x_k,n}(X)|}
%t \epsilon \delta(\|z\|)\gamma(t \|\xi\|) \bar \xi_n,  \\
%\xi_1, \dots, \xi_n \Bigr)
%\end{multline*}
if $k=1,\dots,d$, and
$$
\Theta_t^k (X, z_1, \dots, z_n,\xi_1, \dots, \xi_n) =
(X, z_1, \dots, z_n,\xi_1 - \delta(\|z\|) \rho(t\|\xi\|),
\xi_2, \dots, \xi_n)
$$
%\begin{multline*}
%(X, z_1, \dots, z_n,\xi_1, \dots, \xi_n) \mapsto \\
%\Bigl(
%X, z_1 - \frac {t \epsilon \delta(\|z\|)\gamma(t \|\xi\|) \overline{X_1^k(z)}}
%{\sqrt{ (X_1^k(z))^2 +\dots + (X_n^k(z))^2}} \bar \xi_1, \dots,
%z_n - \frac {t \epsilon \delta(\|z\|)\gamma(t \|\xi\|) \overline{X_n^k(z)}}
%{\sqrt{ (X_1^k(z))^2 +\dots + (X_n^k(z))^2}} \bar \xi_n,  \\
%\xi_1, \dots, \xi_n \Bigr)
%\end{multline*}
if $k=d+1,\dots,d'$.

That is we shift
\begin{equation}   \label{zshift}
\begin{cases}
(z_1,\dots,z_n) \text{ by a vector}  \\
\qquad - t \epsilon \delta(\|z\|)\gamma(t \|\xi\|)
\Bigl( \frac {\overline{\beta_{x_k,1}(X)}}{|\beta_{x_k,1}(X)|} \bar \xi_1,
\dots,
\frac {\overline{\beta_{x_k,n}(X)}}{|\beta_{x_k,n}(X)|} \bar \xi_n \Bigr) &
\text{if $k=1,\dots,d$};  \\
\quad \\
\xi_1 \text{ by a scalar} \quad
- \delta(\|z\|) \rho(t\|\xi\|) &
\text{if $k=d+1,\dots,d'$}.
\end{cases}
\end{equation}

This choice of coefficients
$- \frac {\overline{\beta_{x_k,l}(X)}}{|\beta_{x_k,l}(X)|}$,
the property $\re (X_1^k(z)) > 0$ and the equation (\ref{exponent}) imply that
\begin{equation}  \label{realpart}
\re((\Theta_t^k)^*\langle X, \mu(\zeta) \rangle)
\le  \re(\langle X, \mu(\zeta) \rangle),
\end{equation}
and the equality occurs if and only if
$\Theta_t^k(X,\zeta)=(X,\zeta)$.

We define $\Theta_t: \Omega \times T^*M \to \Omega \times T^*M$ by
$$
\Theta_t= \Theta_t^{d'} \circ \dots \circ \Theta_t^1.
$$
Observe that $\Theta_0$ is the identity map.
The following five lemmas and the proposition are some of the key
properties of $\Theta_t$ that we will use.
We do not give their complete proofs (they can be found in \cite{L1}),
but rather show the key steps only.

\separate

\begin{lem}   \label{equal}
For each $k=1,\dots,d$, the maps $\Theta_t$ and $\Theta_t^k$ coincide
on the set $\{\tilde \psi_k(X,z,\xi);\: X \in \Omega,\: \|z\| < \epsilon\}
\subset \Omega \times T^*M$.
\end{lem}

\pf Follows immediately from condition (\ref{one}).
\qed

\begin{lem} \label{realpartlemma}
If $t>0$ and $\zeta \in T^*M$ does not lie in the zero section,
$$
\re \bigl( (\Theta_t)^*\langle X, \mu(\zeta) \rangle \bigr)
<  \re(\langle X, \mu(\zeta) \rangle).
$$
\end{lem}

\pf By (\ref{realpart}), we have
$$
\re \bigl( (\Theta_t)^*\langle X, \mu(\zeta) \rangle \bigr)
\le  \re(\langle X, \mu(\zeta) \rangle),
$$
and the equality is possible only if $\Theta^k_t(X,\zeta)=(X,\zeta)$
for all $k=1,\dots,d'$. In presence of the condition (\ref{D}) it means that
the equality is possible only if $t=0$ or $\zeta$ lies in the zero section.
\qed

Fix a norm $\|.\|_{T^*M}$ on the cotangent space $T^*M$.

\begin{lem} \label{scaling}
There exists an $R_0>0$ (depending on $t$)
such that whenever $X \in \operatorname{supp}(\phi)$,
$\zeta \in T^*M$ and $\|\zeta\|_{T^*M} \ge R_0$ we have
$\Theta_t (X,E \zeta)=E \Theta_t (X,\zeta)$ for all real $E \ge 1$.
That is $\Theta_t$ almost commutes with scaling the fiber.

Moreover, there is an $\tilde R_0>0$, independent of $t \in (0,1]$,
such that $R_0$ can be chosen to be $\tilde R_0 /t$.
\end{lem}

\pf
Same as the proof of Lemma 18 in \cite{L1}.
Recall that $\Theta_t= \Theta_t^{d'} \circ \dots \circ \Theta_t^1$,
hence it is sufficient to show by induction on $k$, $1 \le k \le d'$,
that there exists an $\tilde R_0>0$ such that whenever
$X \in \operatorname{supp}(\phi)$,
$\zeta \in T^*M$ and $\|\zeta\|_{T^*M} \ge \tilde R_0 /t$,
$$
(\Theta_t^k \circ \dots \circ \Theta_t^1) (X,E \zeta)
=E (\Theta_t^k \circ \dots \circ \Theta_t^1) (X,\zeta)
$$
for all real $E \ge 1$.

Suppose first that $1 \le k \le d$.
When $\|\xi\| > 2/t$, $\gamma (t\|\xi\|) = \frac 1{t\|\xi\|}$ and the shift
vector (\ref{zshift})
\begin{multline*}
- t \epsilon \delta(\|z\|)\gamma(t \|\xi\|)
\Bigl( \frac {\overline{\beta_{x_k,1}(X)}}{|\beta_{x_k,1}(X)|} \bar \xi_1,
\dots,
\frac {\overline{\beta_{x_k,n}(X)}}{|\beta_{x_k,n}(X)|} \bar \xi_n \Bigr)  \\
=
- \frac{\epsilon \delta(\|z\|)}{\|\xi\|}
\Bigl( \frac {\overline{\beta_{x_k,1}(X)}}{|\beta_{x_k,1}(X)|} \bar \xi_1,
\dots,
\frac {\overline{\beta_{x_k,n}(X)}}{|\beta_{x_k,n}(X)|} \bar \xi_n \Bigr)
\end{multline*}
stays unchanged if we replace $(\xi_1,\dots,\xi_n)$ with 
$(E\xi_1,\dots,E\xi_n)$, for any real $E \ge 1$.
Hence in this situation $\Theta_t^k (X,E \zeta)=E \Theta_t^k (X,\zeta)$.

Now suppose that $d < k \le d'$.
When $\|\xi\| > 2/t$, $\rho(t\|\xi\|) = at\|\xi\|$ and the $\xi_1$ coordinate
is shifted by
$$
- \delta(\|z\|) \rho(t\|\xi\|) = - at \delta(\|z\|) \|\xi\|.
$$
It follows that $\Theta_t^k (X,E \zeta)=E \Theta_t^k (X,\zeta)$
whenever $\|\xi\| > 2/t$ and $E \ge 1$.

Set $(X,\zeta_k)= (\Theta_t^{k-1} \circ \dots \circ \Theta_t^1)(X,\zeta)$.
Then one argues by induction on $k$ that there exists an $\tilde R_0>0$
such that whenever $X \in \operatorname{supp}(\phi)$ and
$\|\zeta\|_{T^*M} \ge \tilde R_0 /t$ we have
$\|\xi(X,\zeta_k)\| > 2/t$ which in turn implies
$$
(\Theta_t^k \circ \dots \circ \Theta_t^1) (X,E \zeta)
= \Theta_t^k (X,E\zeta_k) = E \Theta_t^k (X,\zeta_k)
=E (\Theta_t^k \circ \dots \circ \Theta_t^1) (X,\zeta).
$$
\qed

\begin{lem}  \label{kappa}
There exist a smooth bounded function $\tilde\kappa (X, v, t)$ defined on
$$
\Omega \times \{\zeta \in T^*M;\, \|\zeta\|_{T^*M}=1 \} \times [0,1]
$$
and a real number $\tilde r_0>0$
such that, whenever $t \|\zeta\|_{T^*M} \le \tilde r_0$,
$$
\langle X, \mu(\zeta) \rangle -
(\Theta_t^d \circ \dots \circ \Theta_t^1)^* \langle X, \mu(\zeta) \rangle
= t \|\zeta\|_{T^*M}^2 \cdot
\tilde\kappa \bigl( X, \frac{\zeta}{\|\zeta\|_{T^*M}}, t \bigr).
$$
Moreover, $\re(\tilde\kappa)$ is positive and bounded away
from zero for $X \in \supp(\phi)$.
\end{lem}

\pf
Same as the proof of Lemma 22 in \cite{L1}.
Write
\begin{multline*}
\langle X, \mu(\zeta) \rangle -
(\Theta_t^d \circ \dots \circ \Theta_t^1)^* \langle X, \mu(\zeta) \rangle
= \bigl( \langle X, \mu(\zeta) \rangle -
\langle X, \mu(\Theta_t^1(X,\zeta)) \rangle \bigr)   \\
+ \dots + \bigl(
\langle X, \mu((\Theta_t^{k-1} \circ\dots\circ \Theta_t^1)(X,\zeta)) \rangle -
\langle X, \mu((\Theta_t^k \circ \dots \circ \Theta_t^1)(X,\zeta)) \rangle
\bigr)  \\
+ \dots + \bigl( \langle X,
\mu((\Theta_t^{d-1} \circ \dots \circ \Theta_t^1)(X,\zeta)) \rangle -
\langle X,\mu((\Theta_t^d \circ \dots \circ \Theta_t^1)(X,\zeta)) \rangle
\bigr).
\end{multline*}
Let $(X,\zeta_k)= (\Theta_t^{k-1} \circ \dots \circ \Theta_t^1)(X,\zeta)$,
$z_j=z_j(X,\zeta_k)$, $\xi_j=\xi_j(X,\zeta_k)$, $1 \le j \le n$,
and suppose for the moment $t\|\xi\| < 1$ so that $\gamma (t\|\xi\|) = 1$.
Then in the coordinate system $\tilde \psi_k$
\begin{multline*}
\langle X, \mu(\zeta_k) \rangle -
\langle X, \mu(\Theta_t^k(X,\zeta_k)) \rangle  \\
=
t \epsilon \delta(\|z\|)\gamma(t\|\xi\|)
\cdot \bigl( |\beta_{x_k,1}(X)| |\xi_1|^2 + \dots
+ |\beta_{x_k,n}(X)| |\xi_n|^2 \bigr) \\
=
t \|\zeta\|_{T^*M}^2 \epsilon \delta(\|z\|)
\frac{|\beta_{x_k,1}(X)| |\xi_1|^2 + \dots
+ |\beta_{x_k,n}(X)| |\xi_n|^2}{\|\zeta\|_{T^*M}^2}.
\end{multline*}
It is clear that
$\epsilon \delta(\|z\|)
\frac{|\beta_{x_k,1}(X)| |\xi_1|^2 + \dots
+ |\beta_{x_k,n}(X)| |\xi_n|^2}{\|\zeta\|_{T^*M}^2}$ is positive.
Since the set
$\supp(\phi) \times \{\zeta \in T^*M;\, \|\zeta\|_{T^*M}=1 \} \times [0,1]$
is compact, this quotient is bounded away from zero on this set.
Then one argues that there is a real number $\tilde r_0>0$
such that $t \|\zeta\|_{T^*M} \le \tilde r_0$ implies $t\|\xi\| < 1$.
\qed

Similarly we have:

\begin{lem}  \label{kappa'}
There exist a smooth bounded function $\tilde\kappa' (X, v, t)$ defined on
$$
\Omega \times \{\zeta \in T^*M;\, \|\zeta\|_{T^*M}=1 \} \times [0,1]
$$
and a real number $\tilde r'_0>0$
such that, whenever $t \|\zeta\|_{T^*M} \le \tilde r'_0$,
$$
\langle X, \mu(\zeta) \rangle -
(\Theta_t^{d'} \circ \dots \circ \Theta_t^{d+1})^*
\langle X, \mu(\zeta) \rangle
= t^2 \|\zeta\|_{T^*M}^2 \cdot
\tilde\kappa' \bigl( X, \frac{\zeta}{\|\zeta\|_{T^*M}}, t \bigr).
$$
Moreover, $\re(\tilde\kappa')$ is positive and bounded away
from zero for $X \in \supp(\phi)$.
\end{lem}

Finally, the most important property of $\Theta_t$ is:

\begin{prop}  \label{slanting}
For any $t \in [0,1]$, we have:
\begin{multline*}
\lim_{R \to \infty}
\int_{\Omega \times (Ch({\cal F}) \cap \{\|\zeta\|_{T^*M} \le R\} )}
\bigl(
e^{\langle X, \mu(\zeta) \rangle} \phi(X) \wedge \alpha(X) \wedge e^{\sigma} \\
- \Theta_t^*(
e^{\langle X, \mu(\zeta) \rangle} \phi(X) \wedge \alpha(X) \wedge e^{\sigma}
) \bigr) =0.
\end{multline*}
\end{prop}

\pf
The proof of Lemma 19 in \cite{L1} applies here because it is based on
the properties of $\Theta_t$ stated in lemmas \ref{realpartlemma},
\ref{scaling}, \ref{kappa}, \ref{kappa'} and not on any other properties.
It is an integration by parts argument similar to the proof of rapid decay
of the Fourier transform $\hat \phi$ in the imaginary directions.

Since the form
$e^{\langle X, \mu(\zeta) \rangle} \phi(X) \wedge \alpha(X) \wedge e^{\sigma}$
is closed, the integral
\begin{multline*}
\int_{\Omega \times (Ch({\cal F}) \cap \{\|\zeta\|_{T^*M} \le R\} )}
\bigl(
e^{\langle X, \mu(\zeta) \rangle} \phi(X) \wedge \alpha(X) \wedge e^{\sigma}
- \Theta_t^*(
e^{\langle X, \mu(\zeta) \rangle} \phi(X) \wedge \alpha(X) \wedge e^{\sigma}
) \bigr)  \\
=\int_{\Omega \times (Ch({\cal F}) \cap \{\|\zeta\|_{T^*M} \le R\} ) -
(\Theta_t)_* \bigl(
\Omega \times (Ch({\cal F}) \cap \{\|\zeta\|_{T^*M} \le R\} ) \bigr)}
e^{\langle X, \mu(\zeta) \rangle} \phi(X) \wedge \alpha(X) \wedge e^{\sigma}
\end{multline*}
is equal to the integral of
$e^{\langle X, \mu(\zeta) \rangle} \phi(X) \wedge \alpha(X) \wedge e^{\sigma}$
over the chain traced by
$(\Theta_{t'})_* \bigl( \Omega \times
\partial (Ch({\cal F}) \cap \{\|\zeta\|_{T^*M} \le R\}) \bigr)$
as $t'$ varies from $0$ to $t$.
We will show that this integral tends to zero as $R \to \infty$.

Since $Ch({\cal F})$ is a cycle in $T^*M$, the chain
$\Omega \times \partial (Ch({\cal F}) \cap \{\|\zeta\|_{T^*M} \le R\})$
is supported inside the set
$\Omega \times \{\zeta \in T^*M;\, \|\zeta\|_{T^*M} =R \}$.
As $R \to \infty$, we can assume that $R>0$.
Then the chain traced by
$(\Theta_{t'})_* \bigl( \Omega \times
\partial (Ch({\cal F}) \cap \{\|\zeta\|_{T^*M} \le R\}) \bigr)$
as $t'$ varies from $0$ to $t$ lies away from the zero section
$\Omega \times T^*_MM$ in $\Omega \times T^*M$.
If we regard $\Theta$ as a map
$\Omega \times T^*M \times [0,1] \to \Omega \times T^*M$,
we get an integral of
$\Theta^* \bigl(
e^{\langle X, \mu(\zeta) \rangle} \phi(X) \wedge \alpha(X) \wedge e^{\sigma}
\bigr)$
over the chain
$\Omega \times \partial (Ch({\cal F}) \cap \{\|\zeta\|_{T^*M} \le R\})
\times [0,t]$.

The idea is to integrate out the $\Omega$ variable and check that
the result decays faster than any negative power of $R$.
Clearly, $\Theta^*(\phi)=\phi$ and Lemma \ref{realpartlemma} says that
$$
\Theta^* \langle X, \mu(\zeta) \rangle =
\langle X, \mu(\zeta) \rangle - \kappa(X,\zeta, t')
$$
for some smooth function $\kappa(X, \zeta, t')$
which has positive real part.
The integral in question can be rewritten as
$$
\int_{
\Omega \times \partial (Ch({\cal F}) \cap \{\|\zeta\|_{T^*M} \le R\})
\times [0,t]}
e^{\langle X, \mu(\zeta) \rangle}
e^{-\kappa(X, \zeta, t')} \phi(X) \wedge
\Theta^* (\alpha(X) \wedge e^{\sigma}).
$$

We pick a system of local coordinates $(z_1,\dots,z_n)$ of $M$
and construct respective local coordinates
$(z_1,\dots,z_n,\xi_1,\dots,\xi_n)$ of $T^*M$.
Suppose that we know that all the partial
derivatives of all orders of $e^{-\kappa(X, \zeta, t')}$
and $\Theta^* (\alpha(X) \wedge e^{\sigma})$
with respect to the $X$ variable can be bounded independently
of $\zeta$ and $t'$ on the set
$\supp(\phi) \times \{\zeta \in T^*M;\, \|\zeta\|_{T^*M}>0\} \times [0,t]$.
Let $y_1, \dots, y_m$ be a system of linear coordinates on
$\g g_{\BB R}$, write
$\mu(\zeta)=\beta_1(\zeta) dy_1 + \dots +\beta_m(\zeta) dy_m$,
then
\begin{multline*}
\int_{\g g_{\BB R}}
e^{\langle X, \mu(\zeta) \rangle}
e^{-\kappa(X, \zeta, t')} \phi(X) \wedge
\Theta^* (\alpha(X) \wedge e^{\sigma})  \\
= - \frac 1{\beta_l(\zeta)}
\int_{\g g_{\BB R}}
e^{\langle X, \mu(\zeta) \rangle}
\frac {\partial}{\partial y_l} \bigl( e^{-\kappa(X,\zeta, t')}
\phi(X) \wedge \Theta^* (\alpha(X) \wedge e^{\sigma}) \bigr),
\end{multline*}
and the last integral can be bounded by a constant multiple of $R^n$.
We can keep performing integration by parts to get the desired
estimate just like for the ordinary Fourier transform.
Thus, after integrating out the $X$-variable, we see that the
integrand indeed decays rapidly in the fiber variable of $T^*M$.
Hence our integral tends to zero as $R \to \infty$.

To show boundedness of the partial derivatives one follows the proof of
Lemma 19 in \cite{L1} which uses lemmas \ref{realpartlemma},
\ref{scaling}, \ref{kappa} and \ref{kappa'}.
\qed

\separate

Recall the Borel-Moore chain $C(X_0)$ described in Proposition \ref{C(X)prop}.
The set $\Omega$ was chosen so that both
$\Omega$ and $\Omega \cap \g t_{\BB R}(X_0)$ are connected.
Hence, for each $X \in \Omega \cap \g t_{\BB R}(X_0)$,
we can choose $C(X)$ equal $C(X_0)$.
Moreover, for each $X \in \Omega$, we can choose $C(X)$ equal
$\omega(X)_*C(X_0)$. These chains $C(X)$, $X \in \Omega$, piece together
into a Borel-Moore chain in $\Omega \times T^*M$
of dimension $(\dim_{\BB R} \g g_{\BB R}+2n+1)$
which appears in each chart $\tilde \psi_k$ as $\Omega \times C(X_0)$,
$$
\partial C = \Omega \times Ch({\cal F}) -
\Omega \times (m_1(X) T^*_{\omega(X) \cdot x_1}M + \dots
+ m_d(X) T^*_{\omega(X) \cdot x_d}M)
$$
and the support of $C$ lies inside
$\{(X,\zeta) \in \Omega \times T^*M;\:
\re( \langle X,\mu(\zeta) \rangle) \le 0\}$.

Take an $R \ge 1$ and restrict all cycles to the set
$\{ (X,\zeta) \in \Omega \times T^*M;\: \|\zeta\|_{T^*M} \le R\}$.
Let $C_{\le R}$ denote the restriction of the cycle $C$, then
it has boundary
\begin{multline*}
\partial C_{\le R} =
\Omega \times (Ch({\cal F}) \cap \{\|\zeta\|_{T^*M} \le R\}) - C'(R) \\
- \Omega \times
\bigl( m_1(X) (T^*_{\omega(X) \cdot x_1}M \cap \{\|\zeta\|_{T^*M} \le R\})
+ \dots +
m_d (T^*_{\omega(X) \cdot x_d}M \cap \{\|\zeta\|_{T^*M} \le R\}) \bigr),
\end{multline*}
where $C'(R)$ is a $(\dim_{\BB R} \g g_{\BB R} +2n)$-chain supported in the set
$$
\{ (X, \zeta) \in \Omega \times T^*M;\: \|\zeta\|_{T^*M}=R,\:
\re( \langle X, \mu(\zeta) \rangle ) \le 0 \}.
$$
Because the chain $C$ is conic, the piece of boundary $C'(R)$
depends on $R$ by an appropriate scaling of $C'(1)$ in the fiber direction.

\begin{lem}
For a fixed $t \in (0,1]$,
$$
\lim_{R \to \infty} \int_{C'(R)} 
\Theta_t^* \bigl(
e^{\langle X, \mu(\zeta) \rangle} \phi(X) \wedge \alpha(X) \wedge e^{\sigma}
\bigr) =0.
$$
\end{lem}

\pf
Same as the proof of Lemma 20 in \cite{L1}.
Integrating the form
$\Theta_t^* \bigl(
e^{\langle X, \mu(\zeta) \rangle} \phi(X) \wedge \alpha(X) \wedge e^{\sigma}
\bigr)$
over the chain $C'(R)$ is equivalent to integrating
$e^{\langle X, \mu(\zeta) \rangle} \phi(X) \wedge \alpha(X) \wedge e^{\sigma}$
over $(\Theta_t)_* C'(R)$.
Let $R_0$ be as in Lemma \ref{scaling},
then, for $R \ge R_0$, the chain
$(\Theta_t)_* C'(R)$ depends on $R$ by scaling
$(\Theta_t)_* C'(R_0)$ in the fiber direction.
By Lemma \ref{realpartlemma}, for every $(X, \zeta)$ lying in the support of
$C'(R)$, the real part of
$\langle X,\mu(\Theta_t(X, \zeta)) \rangle$ is strictly negative.
By compactness of
$|(\Theta_t)_*C'(R_0)| \cap (\supp(\phi) \times T^*M)$,
there exists an $\epsilon'>0$ such that, whenever $(X, \zeta)$
lies in the support of $(\Theta_t)_*C'(R_0)$ and $X$ lies
in the support of $\phi$, we have
$\re( \langle X,\mu(\zeta) \rangle ) \le -\epsilon'$.
Then, for all $R \ge R_0$ and all
$(X, \zeta) \in
|(\Theta_t)_*C'(R)| \cap (\supp(\phi) \times T^*X)$, we have
$\re( \langle X,\mu(\zeta) \rangle ) \le -\epsilon' \frac R{R_0}$.
Since the integrand decays exponentially over the
support of $(\Theta_t)_* C'(R)$, the integral
tends to zero as $R \to \infty$.
\qed

Thus, using Proposition \ref{slanting},
\begin{multline*}
\int_{Ch({\cal F})}
\mu^* (\widehat{\phi\alpha}) \wedge e^{\sigma}
= \lim_{R \to \infty}
\int_{\Omega \times (Ch({\cal F}) \cap \{\|\zeta\|_{T^*M} \le R\})}
e^{\langle X, \mu(\zeta) \rangle} \phi(X) \wedge \alpha(X) \wedge e^{\sigma} \\
= \lim_{R \to \infty}
\int_{\Omega \times (Ch({\cal F}) \cap \{\|\zeta\|_{T^*M} \le R\})}
\Theta_t^* \bigl(
e^{\langle X, \mu(\zeta) \rangle} \phi(X) \wedge \alpha(X) \wedge e^{\sigma}
\bigr)  \\
=\lim_{R \to \infty}
\int_{C'(R) + \Omega \times \bigl(
{\Sigma}_{k=1}^d m_k(X) (T^*_{\omega(X) \cdot x_k}M \cap
\{\|\zeta\|_{T^*M} \le R\}) \bigr)}
\Theta_t^* \bigl(
e^{\langle X, \mu(\zeta) \rangle} \phi(X) \wedge \alpha(X) \wedge e^{\sigma}
\bigr)  \\
=\lim_{R \to \infty}
\int_{\Omega \times \bigl(
{\Sigma}_{k=1}^d m_k(X) (T^*_{\omega(X) \cdot x_k}M \cap
\{\|\zeta\|_{T^*M} \le R\}) \bigr)}
\Theta_t^* \bigl(
e^{\langle X, \mu(\zeta) \rangle} \phi(X) \wedge \alpha(X) \wedge e^{\sigma}
\bigr),
\end{multline*}
i.e. the integral over $C'(R)$ can be ignored and
we are left with integrals over
$m_k(X) \bigl( \Omega \times
(T^*_{\omega(X) \cdot x_k}M \cap \{\|\zeta\|_{T^*M} \le R\}) \bigr)$,
for $k=1,\dots,d$.
Because the integral converges absolutely, we can let $R \to \infty$ and
drop the restriction $\|\zeta\|_{T^*M} \le R$:
\begin{equation}  \label{integral2}
\int_{Ch({\cal F})}
\mu^* (\widehat{\phi\alpha}) \wedge e^{\sigma} =
\int_{\Omega \times \bigl(
{\Sigma}_{k=1}^d m_k(X) T^*_{\omega(X) \cdot x_k}M \bigr)}
\Theta_t^* \bigl(
e^{\langle X, \mu(\zeta) \rangle} \phi(X) \wedge \alpha(X) \wedge e^{\sigma}
\bigr).
\end{equation}

Lemma \ref{equal} tells us that the maps
$\Theta_t$ and $\Theta_t^k$ coincide over $T^*_{\omega(X) \cdot x_k}M$:
$$
\Theta_t|_{T^*_{\omega(X) \cdot x_k}M} \equiv
\Theta_t^k|_{T^*_{\omega(X) \cdot x_k}M}.
$$
We also have $\delta(\|z\|)=1$, and the exponential part
$\Theta_t^* \bigl( \langle X, \mu(\zeta) \rangle \bigl)$
of our integrand
$$
\Theta_t^* \bigl( e^{\langle X, \mu(\zeta) \rangle}
\phi(X) \wedge \alpha(X) \wedge e^{\sigma} \bigr)
$$
becomes
\begin{equation}  \label{exp}
- t \epsilon \gamma(t\|\xi\|) \bigl( |\beta_1(X)| \xi_1\bar \xi_1 +\dots+
|\beta_n(X)| \xi_n \bar \xi_n \bigr).
\end{equation}

We know that
$\int_{\g g_{\BB R}' \times Ch({\cal F})}
\Theta_t^* \bigl(  e^{\langle X, \mu(\zeta) \rangle}
\phi(X) \wedge \alpha(X) \wedge e^{\sigma} \bigr)$
does not depend on $t$. So in order to calculate its value
we are allowed to regard it as a constant function of $t$
and take its limit as $t \to 0^+$.

We can break up our chain $m_k(X)(T^*_{\omega(X) \cdot x_k}M)$ into
two portions: one portion where $\|\xi(X,\zeta)\| \ge 1/t$ and the
other where $\|\xi(X,\zeta)\| < 1/t$.

\begin{lem}
$$
\lim_{t \to 0^+}
\int_{m_k(X) \bigl( \Omega \times
(T^*_{\omega(X) \cdot x_k}M \cap \{\|\xi(X,\zeta)\| \ge 1/t\}) \bigr)}
\Theta_t^* \bigl( e^{\langle X, \mu(\zeta) \rangle}
\phi(X) \wedge \alpha(X) \wedge e^{\sigma} \bigr) =0.
$$
\end{lem}

\pf
When $\|\xi\| \ge 1/t$,
$\gamma(t \|\xi\|) \ge \frac 1{t \|\xi\|}$
and the exponential part (\ref{exp}) is at most
$$
- \frac {\epsilon}{\|\xi\|} \bigl( |\beta_{x_k,1}(X)| \xi_1\bar \xi_1
+\dots+ |\beta_{x_k,n}(X)| \xi_n \bar \xi_n \bigr).
$$
But $\xi_1 \bar \xi_1 + \dots + \xi_n \bar \xi_n = \|\xi\|^2$, so
at least one of the $\xi_l \bar \xi_l \ge \|\xi\|^2/n$.
Thus we get a new estimate of (\ref{exp}) from above:
$$
- \frac {\epsilon}n |\beta_{x_k,l}(X)| \|\xi\|
\le - \frac {\epsilon}{nt} |\beta_{x_k,l}(X)|.
$$
The last expression tends to $-\infty$ as $t \to 0^+$,
i.e. the integrand decays exponentially and the lemma follows.
\qed

Thus, in the formula (\ref{integral2}) the integral over the portion
$$
m_k(X) \bigl( \Omega \times (T^*_{\omega(X) \cdot x_k}M \cap
\{\|\xi(X,\zeta)\| \ge 1/t\}) \bigr)
$$
can be ignored too:
\begin{multline*}
\int_{Ch({\cal F})}
\mu^* (\widehat{\phi\alpha}) \wedge e^{\sigma}  \\
=\lim_{t \to 0^+}
\int_{ {\Sigma}_{k=1}^d m_k(X) \bigl( \Omega \times
(T^*_{\omega(X) \cdot x_k}M \cap \{\|\xi(X,\zeta)\| < 1/t\}) \bigr)}
\Theta_t^* \bigl(
e^{\langle X, \mu(\zeta) \rangle} \phi(X) \wedge \alpha(X) \wedge e^{\sigma}
\bigr).
\end{multline*}

Finally, over the portion
$m_k(X) \bigl( \Omega \times (T^*_{\omega(X) \cdot x_k}M \cap
\{\|\xi(X,\zeta)\| < 1/t\}) \bigr)$,
the function $\gamma(t \|\xi\|)$ is identically one,
so the exponential part (\ref{exp}) reduces to
$$
- t \epsilon \bigl( |\beta_{x_k,1}(X)| \xi_1\bar \xi_1 +\dots+
|\beta_{x_k,n}(X)| \xi_n \bar \xi_n \bigr).
$$
We also have $\Theta_t^*(\phi)=\phi$,  $\Theta_t^*(d\xi_l)=d\xi_l$,
$$
\Theta_t^*(dz_l) = -d \Bigl( t \epsilon \gamma(t\|\xi\|)
\frac {\overline{\beta_{x_k,l}(X)}} {|\beta_{x_k,l}(X)|} \bar \xi_l \Bigr)
=-t\epsilon \frac {\overline{\beta_{x_k,l}(X)}} {|\beta_{x_k,l}(X)|}
d\bar \xi_l,
\qquad
\Theta_t^*(d \bar z_l) = 
- t\epsilon \frac {\beta_{x_k,l}(X)}{|\beta_{x_k,l}(X)|} d\xi_l,
$$
$$
\Theta_t^*(\sigma) =
- t\epsilon \frac {\overline{\beta_{x_k,1}(X)}} {|\beta_{x_k,1}(X)|}
d \xi_1 \wedge d\bar \xi_1 -\dots-
t\epsilon \frac {\overline{\beta_{x_k,n}(X)}} {|\beta_{x_k,n}(X)|}
d \xi_n \wedge d\bar \xi_n.
$$
The form
$$
\talpha(X)_{[2n]} = 
\Bigl( e^{\langle X, \mu(\zeta) \rangle + \sigma}
\wedge \pi^* \bigl( \alpha(X) \bigr) \Bigr)_{[2n]}
= e^{\langle X, \mu(\zeta) \rangle}
\sum_{l=0}^n \frac 1{l!} \sigma^l \wedge \alpha(X)_{[2n-2l]},
$$
and we end up integrating
\begin{multline*}
e^{-t \epsilon \bigl( |\beta_{x_k,1}(X)| \xi_1\bar \xi_1 +\dots+
|\beta_{x_k,n}(X)| \xi_n \bar \xi_n \bigr)}   \\
\cdot \phi(X) \wedge
\Bigl( (-t \epsilon)^n \Theta_t^* (\pi^* \alpha(X)_{[0]})
\frac {\overline{\beta_{x_k,1}(X)}} {|\beta_{x_k,1}(X)|}
\dots \frac {\overline{\beta_{x_k,n}(X)}}{|\beta_{x_k,n}(X)|}
d\xi_1 \wedge d\bar\xi_1 \wedge \dots \wedge d\xi_n \wedge d\bar\xi_n   \\
+ \text{terms containing
$\Theta_t^* \bigl(\pi^* \alpha(X)_{[2l]} \bigr)$, $l>0$} \Bigr)
\end{multline*}
over $m_k(X) \bigl( \Omega \times (T^*_{\omega(X) \cdot x_k}M \cap
\{\|\xi(X,\zeta)\| < 1/t\}) \bigr)$.
(Recall that the orientation of this chain is determined by the
product orientation on $\Omega \times T^*_{\omega(X) \cdot x_k}M$,
and the orientation of $T^*_{\omega(X) \cdot x_k}M$ is given by
(\ref{orientation}).)

We can write
$$
\Theta_t^* \bigl( \pi^* \alpha(X)_{[0]} \bigr) =
\alpha(X)_{[0]}(\omega(X)\cdot x_k) + 
t \sum_{a=1}^n \bigl( \xi_a A_a(X,t\xi_1,\dots,t\xi_n) +
\bar\xi_a B_a(X,t\xi_1,\dots,t\xi_n) \bigr)
$$
for some bounded functions $A_a$,
$B_a$ of $(X,t\xi_1,\dots,t\xi_n)$, $a=1,\dots,n$.
We can also write
\begin{multline*}
\Theta_t^* \bigl( \pi^* \alpha(X)_{[2]} \bigr)   \\
= t^2 \sum_{b,c=1}^n \bigl( C_{b,c}(X,t\xi_1,\dots,t\xi_n) d\xi_b \wedge d\xi_c
+ D_{b,c}(X,t\xi_1,\dots,t\xi_n) d \bar\xi_b \wedge d\xi_c   \\
+ E_{b,c}(X,t\xi_1,\dots,t\xi_n) d\xi_b \wedge d\bar\xi_c +
F_{b,c}(X,t\xi_1,\dots,t\xi_n) d\bar\xi_b \wedge d\bar\xi_c \bigr),
\end{multline*}
where each of $C_{b,c}$, $D_{b,c}$, $E_{b,c}$, $F_{b,c}$ is a bounded
function in terms of the variables $(X,t\xi_1,\dots,t\xi_n)$.
Similarly we can express
$\Theta_t^* \bigl( \pi^* \alpha(X)_{[2l]} \bigr)$
for $l=1,\dots,n$.
Then, changing variables $y_l = \sqrt{\epsilon t} \xi_l$ for $l=1,\dots,n$,
we obtain the following estimate to (\ref{theintegral}):
\begin{multline*}
(-1)^n\int_{\Omega} m_k(X) \phi(X)
\int_{\{|y_1|^2+\dots+|y_n|^2 < \frac {\epsilon}t \}}
e^{-|\beta_{x_k,1}(X)| |y_1|^2 -\dots- |\beta_{x_k,n}(X)| |y_n|^2}    \\
\cdot \Bigl( \alpha(X)_{[0]}(\omega(X) \cdot x_k)
\frac {\overline{\beta_{x_k,1}(X)}}{|\beta_{x_k,1}(X)|} \dots
\frac {\overline{\beta_{x_k,n}(X)}}{|\beta_{x_k,n}(X)|}
dy_1 \wedge d\bar y_1 \wedge \dots \wedge dy_n \wedge d\bar y_n  \\
+ \sqrt{t} \cdot (\text{bounded terms}) \Bigr).
\end{multline*}
By the Lebesgue dominant convergence theorem this integral tends to
\begin{multline*}
(-1)^n \int_{\Omega} m_k(X) \phi(X) \alpha(X)_{[0]}(\omega(X) \cdot x_k)
\int_{\{(y_1, \dots, y_n) \in \BB C^n\}}
e^{-|\beta_{x_k,1}(X)| |y_1|^2 -\dots- |\beta_{x_k,n}(X)| |y_n|^2}    \\
\cdot \frac {\overline{\beta_{x_k,1}(X)}}{|\beta_{x_k,1}(X)|} \dots
\frac {\overline{\beta_{x_k,n}(X)}}{|\beta_{x_k,n}(X)|}
dy_1 \wedge d\bar y_1 \wedge \dots \wedge dy_n \wedge d\bar y_n  \\
= (-2\pi i)^n \int_{\Omega} m_k(X)
\frac {\alpha(X)_{[0]}(\omega(X) \cdot x_k)}
{\beta_{x_k,1}(X) \dots \beta_{x_k,n}(X)} \phi(X)
\end{multline*}
as $t \to 0^+$.
The last expression may appear to have an extra factor of $(-1)^n$,
but it is correct because of the convention explained in
remarks \ref{T*M} and \ref{orientationremark}.
This proves formula (\ref{mainequation}) when the form $\phi$ is supported
inside $\Omega$.
Then a simple partition of unity argument proves formula (\ref{mainequation})
when the form $\phi$ is compactly supported in $\g g_{\BB R}'$,
an open subset of the set of regular semisimple elements in
$\g g_{\BB R}$ whose complement has measure zero.
Since $\Lambda$ is $G_{\BB R}$-invariant and the form $\alpha$ is
$U_{\BB R}$-equivariant, $F_{\alpha}$ must be invariant under the adjoint
action of $G_{\BB R} \cap U_{\BB R}$.

\separate

To prove the last statement of Theorem \ref{main} we assume that
$F_\alpha$ is a locally $L^1$ function on $\g g_{\BB R}$ and
drop the assumption that the support of $\phi$ lies inside $\g g_{\BB R}'$.
Let $\{\phi_l\}_{l=1}^{\infty}$ be a partition of unity
on $\g g_{\BB R}'$ subordinate to the covering by those open sets
$\Omega$'s.
Then $\phi$ can be realized on $\g g_{\BB R}'$
as a pointwise convergent series:
$$
\phi= \sum_{l=1}^{\infty} \phi_l \phi.
$$
Because $F_{\alpha} \in L^1_{loc} (\g g_{\BB R})$, the series
$\sum_{l=1}^{\infty}\int_{\g g_{\BB R}} F_{\alpha} \phi_l \phi$
converges absolutely. Hence
\begin{multline*}
\int_{Ch({\cal F})} \mu^* (\widehat{\phi\alpha}) \wedge e^{\sigma} =
\int_{Ch({\cal F})} \Bigl( \int_{\g g_{\BB R}} \talpha \wedge \phi(X) \Bigr) \\
= \sum_{l=1}^{\infty} \int_{Ch({\cal F})}
\Bigl( \int_{\g g_{\BB R}} \talpha \wedge \phi_l \phi(X) \Bigr)
= \sum_{l=1}^{\infty}\int_{\g g_{\BB R}} F_{\alpha} \phi_l \phi =
\int_{\g g_{\BB R}} F_{\alpha} \phi,
\end{multline*}
which completes our proof of Theorem \ref{main}.
\qed

\separate

\begin{section}
{A Gauss-Bonnet Theorem for Constructible Sheaves}
\end{section}

In this section we use Theorem \ref{main} to prove a generalization
of the Gauss-Bonnet Theorem for constructible sheaves.

As before, let $G_{\BB C}$ be a connected complex algebraic reductive group
which is defined over $\BB R$, and let $G_{\BB R}$ be a subgroup of
$G_{\BB C}$ lying between the group of real points $G_{\BB C}(\BB R)$
and the identity component $G_{\BB C}(\BB R)^0$.
Let $\g g_{\BB C}$ and $\g g_{\BB R}$ be their respective Lie algebras.
This time we require $U_{\BB R} \subset G_{\BB C}$ to be a compact real
form of $G_{\BB C}$, and let $\g u_{\BB R}$ denote its Lie algebra.
As before, $M$ is a smooth complex projective variety with a complex
algebraic $G_{\BB C}$-action on it such that a maximal complex torus
$T_{\BB C} \subset G_{\BB C}$ acts on $M$ with isolated fixed points,
and ${\cal F}$ is a $G_{\BB R}$-equivariant sheaf on $M$ with
$\BB R$-constructible cohomology.
We assume that the holomorphic moment map
$\mu: T^*M \to \g g_{\BB C}^*$ is proper on the set
$\supp(\sigma|_{Ch({\cal F})})$.
Let $n = \dim_{\BB C} M$.

Pick a $U_{\BB R}$-invariant connection $\nabla$ on the tangent bundle $TM$.
Then N.~Berline, E.~Getzler and M.~Vergne define in Section 7.1 of \cite{BGV}
the equivariant connection and the
{\em equivariant curvature $F_{\g u_{\BB R}}$} associated to $\nabla$.
After that they define the {\em equivariant Euler form}
$$
\chi_{\g u_{\BB R}} (\nabla)(X) = {\det}^{1/2} (- F_{\g u_{\BB R}} (X)),
\qquad X \in {\g u_{\BB R}}.
$$
The form $\chi_{\g u_{\BB R}} (\nabla)$ is $U_{\BB R}$-equivariantly closed
and its class in equivariant cohomology does not depend on the choice of the
$U_{\BB R}$-invariant connection $\nabla$.
It is easy to see that the map
$\chi_{\g u_{\BB R}} (\nabla): \g u_{\BB R} \to \Omega^*(M)$
is polynomial and extends uniquely to a holomorphic polynomial
(but not $G_{\BB C}$-equivariant) function
$$
\chi_{\g g_{\BB C}}:
\g g_{\BB C} \simeq \g u_{\BB R} \otimes_{\BB R} \BB C \to \Omega^*(M).
$$
We use the following properties of $\chi_{\g g_{\BB C}}$:
$$
\chi_{\g g_{\BB C}} (X)_{[2n]} = \text{Euler form of $TM$},
\qquad \forall X \in \g g_{\BB C};
$$
$$
\chi_{\g g_{\BB C}} (X)_{[2k]} \in \Omega^{(k,k)} (M),
\qquad \forall k \in \BB N;
$$
if $p \in M_0(X)$, then
\begin{equation}  \label{chi}
\chi_{\g g_{\BB C}} (X)_{[0]} (p) = i^n \cdot \den_p(X);
\end{equation}
in particular, $\chi_{\g g_{\BB C}}$ satisfies the Conditions \ref{conditions}.

\separate

\begin{thm}  \label{bonnet}
Under the above conditions, if $\phi$ is a smooth compactly supported
differential form on $\g g_{\BB R}$ of top degree,
$$
(2\pi)^{-\dim_{\BB C} M} \int_{Ch({\cal F})}
\mu^* (\widehat{\phi \, \chi_{\g g_{\BB C}}}) \wedge e^{\sigma}
= (2\pi)^{-\dim_{\BB C} M} \int_{Ch({\cal F})} \Bigl( \int_{\g g_{\BB R}}
\widetilde{\chi_{\g g_{\BB C}}} \wedge \phi(X) \Bigr)
= \chi(M, {\cal F}) \cdot \int_{\g g_{\BB R}} \phi,
$$
where
$$
\widetilde{\chi_{\g g_{\BB C}}} (X) = 
e^{\langle X, \mu(\zeta) \rangle + \sigma}  \wedge \chi_{\g g_{\BB C}}(X),
$$
$\chi(M, {\cal F})$ is the Euler characteristic of $M$ with respect to
${\cal F} \in \C (M)$.
\end{thm}

\begin{rem}
If ${\cal F}$ is the constant sheaf $\BB C_M$ on $M$, then
$Ch({\cal F}) = [M]$, the moment map $\mu$ is automatically proper on
$|Ch({\cal F})| =M$, and we obtain the classical Gauss-Bonnet theorem
$$
\chi(M) = (2\pi)^{-\frac 12 \dim_{\BB R} M} \int_M \text{ Euler class of $TM$}.
$$
Here we do not even need the requirement that a maximal complex torus
$T_{\BB C} \subset G_{\BB C}$ acts on $M$ with isolated fixed points.
\end{rem}

\pf
First, we assume that the support of the test form $\phi$ lies in
$\g g_{\BB R}'$.
An immediate application of Theorem \ref{main} together with
the property (\ref{chi}) show that
$$
(2\pi)^{-\dim_{\BB C} M} \int_{Ch({\cal F})}
\mu^* (\widehat{\phi \, \chi_{\g g_{\BB C}}}) \wedge e^{\sigma}
= (2\pi)^{-\dim_{\BB C} M} \int_{Ch({\cal F})} \Bigl( \int_{\g g_{\BB R}}
\widetilde{\chi_{\g g_{\BB C}}} \wedge \phi(X) \Bigr)
= \int_{\g g_{\BB R}} E \phi,
$$
where, using the global coefficient formula (\ref{m2}),
$$
E(X) = \sum_{x_k \in M_0(X)} m_k(X) =
\sum_{x_k \in M_0(X)} \chi(M, {\cal F}_{O_k})
= \chi(M, {\cal F}).
$$

Finally, the constant function $\chi(M, {\cal F})$ is clearly locally
integrable with respect to the Lebesgue measure on $\g g_{\BB R}$,
hence the last part of Theorem \ref{main} applies here and this proves
Theorem \ref{bonnet} in general.
\qed

\separate

\begin{section}
{Duistermaat-Heckman Measures}
\end{section}

%As before, $M$ is a smooth complex projective variety with a complex
%algebraic $G_{\BB C}$-action on it such that a maximal complex torus
%$T_{\BB C} \subset G_{\BB C}$ acts on $M$ with isolated fixed points,
%and $\Lambda \in {\cal L}^+_{G_{\BB R}}(M)$ is a conic Lagrangian
%$G_{\BB R}$-invariant cycle in $T^*M$.

As before, $G_{\BB R}$ is a linear real reductive Lie group with
complexification $G_{\BB C}$,
we denote by $\g g_{\BB R}$ and $\g g_{\BB C}$ their respective
Lie algebras.
We pick another subgroup $U_{\BB R}$ of $G_{\BB C}$ such that, letting
$\g u_{\BB R}$ be the Lie algebra of $U$, we have an isomorphism
$\g u_{\BB R} \otimes_{\BB R} \BB C \simeq \g g_{\BB C}$.
For instance, $U_{\BB R}$ may equal $G_{\BB R}$, but in most interesting
situations $U_{\BB R}$ is a compact real form of $G_{\BB C}$.

Let $M$ be a smooth complex projective variety equipped with an algebraic
action of $G_{\BB C}$ preserving a complex-valued 2-form $\omega$,
and suppose that the restriction of the
$G_{\BB C}$-action to $U_{\BB R}$ is Hamiltonian with respect to $\omega$.
In other words, there exists a moment map
$J: M \to \g u_{\BB R}^* \otimes_{\BB R} \BB C \simeq \g g_{\BB C}^*$
such that
$$
\iota(X_M) \omega = dJ(X),
\qquad \forall X \in \g u_{\BB R}.
$$
Note that we do not require the 2-form $\omega$ to be symplectic,
i.e. $\omega^{\dim_{\BB R} M/2} \ne 0$.
Even the case $\omega=0$, $J=0$ is interesting enough,
but, of course, symplectic forms are the most interesting ones.
We can regard $J: M \to \g g_{\BB C}^*$ as a map
$J: \g g_{\BB C} \to {\cal C}^{\infty}(M)$.
Then $\omega +J$ is an equivariantly closed form on $M$
for the action of $U_{\BB R}$.

Recall that $\sigma$ denotes the canonical complex algebraic holomorphic
symplectic form on the holomorphic cotangent bundle $T^*M$ and
$\mu: T^*M \to \g g_{\BB C}^*$ is the ordinary holomorphic moment map.
Let $\Lambda \in {\cal L}^+_{G_{\BB R}}(M)$ be a conic real-Lagrangian
$G_{\BB R}$-invariant cycle in $T^*M$.
As before, $n = \dim_{\BB C} M$. The Liouville form
$$
\frac {(\omega+\sigma)^n}{n!} =
\bigl( \exp (\omega+\sigma) \bigr)_{[2n]}
$$
determines a measure $\beta_{\Lambda}$ on $\Lambda$.
We call the pushforward of this measure
$(J+\mu)_*(\beta_{\Lambda})$ on $\g g_{\BB C}^*$
the Duistermaat-Heckman measure.
That is, for a compactly supported smooth function
$f \in {\cal C}^{\infty}_c (\g g_{\BB C}^*)$,
\begin{equation}  \label{measure}
\int_{\g g_{\BB C}^*} f \, d(J+\mu)_*(\beta_{\Lambda}) \quad =_{def}
\quad \int_{\Lambda} \frac {(\omega+\sigma)^n}{n!}
\bigl( f \circ (J+\mu) \bigr).
\end{equation}
The right hand side of (\ref{measure}) converges whenever the map
$J+\mu$ is proper on the set $\supp(\sigma|_{\Lambda})$.
This happens whenever $\mu$ is proper on $\supp(\sigma|_{\Lambda})$.
In particular, the pushforward $(J+\mu)_*(\beta_{\Lambda})$
is well-defined when $\mu$ is proper on $|\Lambda|$.

Duistermaat-Heckman measures are important invariants of symplectic manifolds
and there are so many papers on this subject that it is impossible to list
them all.
At first an explicit formula was given by J.~J.~Duistermaat and G.~J.~Heckman
\cite{DH} using the method of exact stationary phase in the special case when
$G$ is a compact torus acting with isolated fixed points.
It was extended to compact non-abelian groups by V.~Guillemin and E.~Prato
\cite{GP}. Then it was extended to compact non-abelian groups acting with
possibly non-isolated fixed points by L.~Jeffrey and F.~Kirwan \cite{JK}.
Many recent results on Duistermaat-Heckman measures are obtained by
computing their Fourier transforms using the integral localization formula
and then inverting these Fourier transforms.

Since the cycle $\Lambda$ is real-Lagrangian and $G_{\BB R}$-invariant,
the moment map $\mu$ takes purely imaginary values on its support $|\Lambda|$:
$$
\mu(|\Lambda|) \quad \subset \quad i\g g_{\BB R}^* \quad \subset
\quad \g g_{\BB R}^* \oplus i \g g_{\BB R}^* \quad \simeq \quad
\g g_{\BB C}^*.
$$
Since $M$ is compact, the support of $(J+\mu)_*(\beta_{\Lambda})$,
which must lie inside $(J+\mu)_*(|\Lambda|)$, is a subset of 
$\g g_{\BB C}^* \simeq \g g_{\BB R}^* \oplus i \g g_{\BB R}^*$
with bounded real part.

The Fourier transform of the Duistermaat-Heckman measure is a distribution
on $\g g_{\BB R}$, i.e. a continuous linear functional on the space
$\Omega^{top}_c(\g g_{\BB R})$ consisting of differential forms of top
degree on $\g g_{\BB R}$ with compact support.
For $\phi \in \Omega^{top}_c(\g g_{\BB R})$, its Fourier transform
$\hat \phi$ is defined by (\ref{ftransform});
recall that $\hat \phi(\xi)$ decays rapidly as $\xi \to \infty$ and
the real part of $\xi$ stays uniformly bounded.
Hence the value of the Fourier transform of $(J+\mu)_*(\beta_{\Lambda})$
at $\phi \in \Omega^{top}_c(\g g_{\BB R})$ is
\begin{multline}  \label{FT1}
\widehat{(J+\mu)_*(\beta_{\Lambda})}(\phi) =
\int_{\g g_{\BB C}^*} \Bigl(
\int_{\g g_{\BB R}} e^{\langle X, \xi \rangle} \phi(X) \Bigr)
\, d(J+\mu)_*(\beta_{\Lambda})   \\
= \int_{\Lambda} \Bigl( \int_{\g g_{\BB R}}
e^{\langle X, (J + \mu)(\zeta) \rangle} \phi(X) \Bigr)
\frac {(\omega+\sigma)^n}{n!},
\qquad X \in \g g_{\BB R}, \: \zeta \in |\Lambda| \subset T^*M.
\end{multline}
We introduce a $U_{\BB R}$-equivariant form
$\alpha: \g g_{\BB C} \to \Omega^*(M)$:
$$
\alpha(X) = \exp(J(X)+ \omega),
$$
then (\ref{FT1}) can be rewritten as
\begin{multline}  \label{FT2}
\widehat{(J+\mu)_*(\beta_{\Lambda})}(\phi) =
\int_{\Lambda} \Bigl( \int_{\g g_{\BB R}} 
e^{\langle X, \mu(\zeta) \rangle + \sigma}
\wedge \phi(X) \wedge \alpha(X) \Bigr)_{[\dim_{\BB R}M]}  \\
= \int_{\Lambda} \Bigl( \int_{\g g_{\BB R}} \talpha \wedge \phi(X) \Bigr),
\qquad X \in \g g ,\: \zeta \in |\Lambda| \subset T^*M.
\end{multline}
This integral is exactly of type (\ref{int}), hence convergent.
The generalized localization formula (\ref{mainequation}) immediately implies:

\begin{prop}  \label{DH}
Suppose there exists a maximal complex torus
$T_{\BB C} \subset G_{\BB C}$ acting on $M$ with finitely many isolated
fixed points and that
$$
\omega \in \Omega^{(2,0)}(M) \oplus \Omega^{(1,1)}(M).
$$
Then the restriction of the Fourier transform of
the Duistermaat-Heckman measure (\ref{FT2}) to $\g g'_{\BB R}$ equals
$$
\widehat{(J+\mu)_*(\beta_{\Lambda})}(\phi) =
\int_{\g g_{\BB R}} F_{\omega}(X) \phi(X),
$$
where $F_{\omega}$ is an $Ad(G_{\BB R} \cap U_{\BB R})$-invariant
function on $\g g'_{\BB R}$ given by the formula
$$
F_{\omega}(X) = (-2\pi)^{\dim_{\BB R} M/2}
\sum_{p \in M_0(X)} m_p(X) \frac {e^{\langle X, J(p) \rangle}}
{\den_p(X)},
$$
where $M_0(X)$ is the set of zeroes of the vector field $X_M$ on $M$,
and $m_p(X)$'s are certain integer multiplicities given by formula (\ref{m2}).
\end{prop}

Note that this formula for $\widehat{(J+\mu)_*(\beta_{\Lambda})}$
is non-trivial even when $\omega=0$, $J=0$.

%Since we want to apply the generalized localization formula
%(\ref{mainequation}), we assume that there exists a maximal complex torus
%$T_{\BB C} \subset G_{\BB C}$ acting on $M$ with finitely many isolated
%fixed points and that
%$$
%\omega \in \Omega^{(2,0)}(M) \oplus \Omega^{(1,1)}(M).
%$$
%Then (\ref{mainequation}) says that the restriction of the Fourier
%transform of the Duistermaat-Heckman measure (\ref{FT2}) to
%$\g g'_{\BB R}$ equals
%$$
%\widehat{(J+\mu)_*(\beta_{\Lambda})}(\phi) =
%\int_{\g g_{\BB R}} F_{\omega}(X) \phi(X),
%$$
%where $F_{\omega}$ is an $Ad(G_{\BB R} \cap U_{\BB R})$-invariant
%function on $\g g'_{\BB R}$ given by the formula
%$$
%F_{\omega}(X) = (-2\pi)^{\dim_{\BB R} M/2}
%\sum_{p \in M_0(X)} m_p(X) \frac {e^{\langle X, J(p) \rangle}}
%{\den_p(X)},
%$$
%where $M_0(X)$ is the set of zeroes of the vector field $X_M$ on $M$,
%and $m_p(X)$'s are certain integer multiplicities given by formula (\ref{m2}).
%Note that this formula is non-trivial even when $\omega=0$, $J=0$.

\separate

\noindent
{\em E-mail address:} {matvei.libine@yale.edu}

\noindent
{\em Department of Mathematics, Yale University,
P.O. Box 208283, New Haven, CT 06520-8283}

\end{document}